\documentclass[reqno,12pt]{amsart}
\usepackage{color}
\usepackage{amsthm}
\usepackage{amsmath}
\usepackage{epsfig}
\usepackage{graphicx}
\usepackage{amssymb}
\usepackage{latexsym}
\usepackage{pdfpages}

\usepackage{enumitem}

\usepackage{ulem} 

\usepackage{ifpdf}
\newcommand{\res}{\!\!\mathop{\hbox{
                                \vrule height 7pt width .5pt depth 0pt
                                \vrule height .5pt width 6pt depth 0pt}}
                                \nolimits}

\def\z{{\bf z}}

\ifpdf 
  \usepackage[hidelinks]{hyperref}
\else 
\fi 

 \usepackage[usenames,dvipsnames]{pstricks}
 \usepackage{epsfig}
 \usepackage{pst-grad} 
 \usepackage{pst-plot} 

\usepackage{color}

\newtheorem{theorem}{Theorem}[section]
\newtheorem{lemma}[theorem]{Lemma}
\newtheorem{definition}[theorem]{Definition}
\newtheorem{proposition}[theorem]{Proposition}
\newtheorem{corollary}[theorem]{Corollary}
\newtheorem{remark}[theorem]{Remark}
\newtheorem{example}[theorem]{Example}
\newtheorem*{theorem*}{\it Theorem}

\def\vint_#1{\mathchoice%
          {\mathop{\kern 0.2em\vrule width 0.6em height 0.69678ex depth -0.58065ex
                  \kern -0.8em \intop}\nolimits_{\kern -0.4em#1}}%
          {\mathop{\kern 0.1em\vrule width 0.5em height 0.69678ex depth -0.60387ex
                  \kern -0.6em \intop}\nolimits_{#1}}%
          {\mathop{\kern 0.1em\vrule width 0.5em height 0.69678ex
              depth -0.60387ex
                  \kern -0.6em \intop}\nolimits_{#1}}%
          {\mathop{\kern 0.1em\vrule width 0.5em height 0.69678ex depth -0.60387ex
                  \kern -0.6em \intop}\nolimits_{#1}}}
\def\vintslides_#1{\mathchoice%
          {\mathop{\kern 0.1em\vrule width 0.5em height 0.697ex depth -0.581ex
                  \kern -0.6em \intop}\nolimits_{\kern -0.4em#1}}%
          {\mathop{\kern 0.1em\vrule width 0.3em height 0.697ex depth -0.604ex
                  \kern -0.4em \intop}\nolimits_{#1}}%
          {\mathop{\kern 0.1em\vrule width 0.3em height 0.697ex depth -0.604ex
                  \kern -0.4em \intop}\nolimits_{#1}}%
          {\mathop{\kern 0.1em\vrule width 0.3em height 0.697ex depth -0.604ex
                  \kern -0.4em \intop}\nolimits_{#1}}}

\def\R{\mathbb R}

\numberwithin{equation}{section}

\def\1{\raisebox{2pt}{\rm{$\chi$}}}

\definecolor{violet(ryb)}{rgb}{0.53, 0.0, 0.69}
\definecolor{internationalorange}{rgb}{1.0, 0.31, 0.0}
\usepackage{mathtools}
\mathtoolsset{showonlyrefs}

\begin{document}

\title[Cahn-Hilliard Equations]{\bf  Cahn-Hilliard Equations  on  Random Walk Spaces}

\author[J. M. Maz\'on and J. Toledo]{Jos\'e M. Maz\'on and Juli\'{a}n Toledo}

 \address{ J. M. Maz\'{o}n \& J. Toledo: Departamento de An\'{a}lisis Matem\'{a}tico,
Univ. Valencia, Dr.~Moliner 50, 46100 Burjassot, Spain;
 {\tt mazon@uv.es},
 {\tt toledojj@uv.es}.}



\setcounter{tocdepth}{1}


\begin{abstract}
In this paper we study a  nonlocal Cahn-Hilliard model  (CHE) in the framework of random walk spaces, which includes as particular cases, the CHE on locally finite weighted connected graphs, the CHE determined by finite Markov chains or the Cahn-Hilliard Equations driven by convolution  integrable kernels.  We consider different transitions for the phase and the chemical potential, and a large class of potentials including obstacle ones. We prove existence and uniqueness of solutions  in $L^1$ of the Cahn-Hilliard Equation. We also show that the Cahn-Hilliard equation is the gradient flow of the Ginzburg-Landau free energy functional on an appropriate Hilbert space. We finally study the asymptotic behaviour of the solutions.

\end{abstract}

\maketitle


{ \renewcommand\contentsname{Contents}
\setcounter{tocdepth}{3}
\addtolength{\parskip}{-0.2cm}
{\small \tableofcontents}
\addtolength{\parskip}{0.2cm} }

\section{Introduction}

  The {\it Cahn-Hilliard model} was formulated by J. W. Cahn and J. E. Hilliard (\cite{CH})   to describe the phase
separation of a binary fluid or alloy. Let $\Omega \subset \R^N$ be the subset where the phase separation takes place. Let $u$ be the concentration of a substance, which takes the values in $[ -1, 1]$. The pure phases of the material correspond to $u=1$ and $u=-1$   while $u \in (-1,1)$ corresponds to the transition in the interface between such phases. Denoting by $\mu$ to the chemical potential, the model, in a simplified version, is described    by the following coupled system of equations
\begin{equation}\label{Ch1}
\left\{   \begin{array}{ll} \displaystyle\frac{\partial u}{\partial t}  -\Delta \mu =0 \quad &\hbox{in} \ (0,\infty) \times \Omega, \\[12pt] \mu = - \epsilon^2 \Delta u + F'(u) & \hbox{in} \ (0,\infty) \times \Omega, \end{array} \right.
\end{equation}
 joint with homogeneous Neumann boundary conditions and initial data, where $F(u)$ is a {\it double-well potential},  where   $\epsilon>0$ is a small  interaction  parameter related to the length of the interface.

 A physically relevant choice for $F$ is a singular logarithmic double-well potential
\begin{equation}\label{logaritmo}F_1(r) = ((1+r) \log (1+r)+ (1-r)  \log (1-r)) - \frac{b}{2}r^2, \quad r \in (-1,1), \ \ \   b >2,
\end{equation}
which is often approximate by regular  double-well polynomial potentials, like  \begin{equation}\label{doblewell}F_2(r) = \frac14(r^2 -1)^2.
\end{equation}
Another physically relevant choice for $F$ is the double-well obstacle potential
\begin{equation}\label{obstacle}
F_3(r):= \frac{a}{2}(1 -r^2) + I_{[-1,1]}(r)= \left\{ \begin{array}{ll}  \displaystyle \frac{a}{2}(1 -r^2) \quad &\hbox{if} \ \vert r \vert \leq 1, \\[8pt] +\infty \quad &\hbox{if} \ \vert r \vert > 1, \end{array} \right.
\end{equation}
where $I_{[-1,1]}$ is the indicator function of $[-1,1]$ and   $a >0$.
 Here, $F_3'(r)$ must be understood as $\partial F_3(r)$, the subdifferential of $F_3$ at $r$.

   Observe that, for the logarithmic potential $F_1$ and the obstacle potential $F_3$, the system  provides a solution $u\in[-1,1]$ in the admissible range, which  does not happen with $F_2$. On the other hand, for the logarithmic,   $u$ can not attain pure phases while this is possible  in the case of the obstacle potential $F_3$.

     As pointed out by   Fife in~\cite{Fife},   the Cahn-Hilliard model \eqref{Ch1} is the $H^{-1}$-gradient flow of the following {\it Ginzburg-Landau} free energy
$$ \mathcal{E} (u) = \frac{\epsilon^2}{2} \int_\Omega \vert \nabla u \vert^2 dx + \int_\Omega F(u) dx.$$
 With the term $  \displaystyle\frac{\epsilon^2}{2} \int_\Omega \vert \nabla u \vert^2 dx$, the model is assuming  a short-range  interaction between  particles and penalties sudden changes for the concentration measure.

 From a microscopic model for lattice gas, Giacomin and Lebowitz \cite{GL}  derived   a nonlocal version of the Cahn-Hilliard system   that takes into account long-interactions between particles.   Simplifying their model, they   consider the nonlocal free energy
\begin{equation}
\mathcal{E}_{NL}(u) = \frac14 \int_\Omega J(x,y)  \vert u(y) - u(x) \vert^2  dxdy + \int_\Omega F(u) dx,
\end{equation}
where  $J(x,y)=\overline J(\vert y - x \vert)$  is a non-singular kernel,
 and associate to it the following  local-nonlocal  Cahn-Hilliard system \begin{equation}\label{ChNL1}
\left\{   \begin{array}{ll}\displaystyle
 \frac{\partial u}{\partial t}  + \Delta \mu =0 \quad &\hbox{in} \ (0,\infty) \times \Omega, \\[14pt]
\displaystyle \mu(x)=-\int_\Omega J(x,y)( u(y)-u(x)) dy + F'(u(x)) &  \hbox{in} \ (0,\infty) \times \Omega, \end{array} \right.
\end{equation}
  joint with local homogeneous Neumann boundary conditions for the first equation (nonlocal regional Neumann boundary conditions are implicit  for the second one) and initial data.

The Cahn-Hilliard model and its variants have been widely
used in different areas of science, see, e.g., \cite{BEG}, \cite{GLSS}, \cite{WLFC} and the monograph \cite{Miranville}.
There is also an extensive literature about nonlocal Cahn-Hilliard equations, see, e.g., \cite{BH}, \cite{BG},  \cite{CoSp}, \cite{DRST}, \cite{DSTARMA}, \cite{DSTJDE},  \cite{GZ}, \cite{GGG}, \cite{GL1}, and the references cited therein.  For an overview of early and recent
references and extensions of the nonlocal Cahn-Hilliard model, we refer  to \cite{BH}, \cite{Miranville0}, \cite{Miranville}.

In \cite{Gal1} Gal considers a double nonlocal Cahn-Hilliard equations in which also the  equation of motion for mass transport  is nonlocal in terms of an operator of the form
$$\mathcal{L}(\mu) = 2 P.V. \int_\Omega K(x,y) (\mu(y) - \mu(x)) dy,$$
with a singular kernel $K(x,y)=\overline K(\vert x - y \vert)$;
 therefore also allowing long-range interaction to occur for $\mu$ between any two points $x,y \in \Omega$.   The governing system of equations   takes the form
 \begin{equation}\label{ChNL2}
\left\{   \begin{array}{ll} \displaystyle\frac{\partial u}{\partial t} = \mathcal{L}(\mu)  \quad &\hbox{in} \ (0,\infty) \times \Omega, \\[12pt]
\displaystyle \mu(x)=-\int_\Omega J(x,y)( u(y)-u(x)) dy + F'(u(x)) & \hbox{in} \ (0,\infty) \times \Omega, \end{array} \right.
\end{equation}
  which is named as a {\it strong-to-weak} interaction   Cahn-Hilliard system.  Recently, Gal and Shomberg  have dealt with a {\it weak-to-weak} interaction Cahn-Hilliard system in \cite{GS}, by considering both $\overline K, \overline J \in L^1(\R^N)$.

Our  aim is to study   Cahn-Hilliard systems in the framework of  random walk  spaces, that has as a particular cases  (double) nonlocal Cahn-Hilliard systems    for integrable kernels, that is, the weak-to-weak iteration, and    Cahn-Hilliard systems in weighted graphs.
  We consider different transitions  for the phase and for the chemical potential.

  Let $[X,\mathcal{B},m^i,\nu_i]$, $i=1,2$, be random walk spaces with $\nu_i$ reversible, and $m_i$-connected (all the concepts used in this introduction can be found in the Preliminaries section).  We are interested in the study of the generalized (doubly nonlocal) Cahn-Hilliard system stated on random walk spaces:
  \begin{equation}\label{DNLS1preforini}
\left\{\begin{array}{lll} \displaystyle u_t(t,x) = \Delta_{m^1}\mu(t,x) ,   &(t,x) \in (0, \infty) \times X, \\[10pt]
\displaystyle \mu(t,x)\in - \Delta_{m^2}u(t,x) + \partial F(u),   &(t,x) \in (0, \infty) \times X,
\\[10pt] u(0,x) = u_0(x),   &x \in  X,\end{array} \right.
\end{equation}
for the (nonlocal) $m_i$-Laplacians ($i=1,2$):
$$\Delta_{m^i}u(t,x)=\int_\Omega  (u(t,y)-u(t,x))dm^i_x(y),$$
and for $F:\mathbb{R}\to ]-\infty,+\infty]$ in the   class of potentials for which its subdifferential (see~\eqref{subdif1}), is
\begin{equation}\label{dondelac}\partial F(r)  = \gamma^{-1}(r) - cr,
\end{equation}
with $\gamma$  any maximal monotone graph in $\R \times  \R$ with $0\in\gamma(0)$  and    $\inf\{\mbox{Ran}(\gamma)\}<\sup\{\mbox{Ran}(\gamma)\}$,  and  with $c>0$.  This class  encompasses a large set of potentials due to the generality of the conditions on $\gamma$. Observe the three choices  of potential given at the beginning belong to such class since:
$$\partial F_1(r)= F_1'(r) = \log(1+r) - \log(1-r) -br,$$
and $\log(1+r) - \log(1-r)$ is the corresponding monotone function;
$$\partial F_2(r)=F'_2(r) = r^3 -r,$$
and $r^3$ is the corresponding  monotone function;
and  $$ \partial F_3 (r) = \partial I_{[-1,1]}(r) - ar,$$
and $\partial I_{[-1,1]}$ is the corresponding  (multivalued) monotone graph
$$\partial I_{[-1,1]}(r) = \left\{ \begin{array}{lll} (-\infty,0] \quad &\hbox{if} \ \ r=-1, \\[5pt] 0  \quad &\hbox{if} \ \ r \in ]-1,1[, \\[5pt] [0, + \infty)  \quad &\hbox{if} \ \ r=1. \end{array} \right. $$

Observe that~\eqref{DNLS1preforini} can be written as follows:
\begin{equation}\label{DNLS1preconX}
\left\{\begin{array}{lll} u_t(t,x) = \Delta_{m^1}  \mu(t,x) ,  &(t,x) \in (0, \infty) \times X, \\[10pt] \mu(t,x) \in - \Delta_{m^2}u(t,x)  +v(t,x)-cu(x,t),   &(t,x) \in (0, \infty) \times X,
\\[14pt]
u(t,x)\in \gamma(v(t,x)), &(t,x) \in (0, \infty) \times X, \\[10pt] u(0,x) = u_0(x), \ &x \in  X,\end{array} \right.
\end{equation}
 that can be reduced  to the study of
\begin{equation}\label{1901}
\left\{\begin{array}{ll} u_t = \Delta_{m^1} v  -\Delta_{m^1\ast m^2} u + (1-c)\Delta_{m^1}u + \Delta_{m^2}u &\hbox{in } (0, \infty) \times X,
\\[10pt]
u \in \gamma(v) &\hbox{in } (0, \infty) \times X,
 \\[10pt] u(0)= u_0 &\hbox{in }  X.\end{array} \right.
\end{equation}
In this last formulation,  the evolution equation
 \begin{equation}\label{Ch1laecnolo}u_t = \Delta_{m^1} v  -\Delta_{m^1\ast m^2} u + (1-c)\Delta_{m^1}u + \Delta_{m^2}u
 \end{equation}
corresponds to the  nonlocal analogous of the Cahn-Hilliard equation
 \begin{equation}\label{Ch1laec}
 \displaystyle u_t  -\Delta( -  \Delta u + F'(u)) =0  \quad \hbox{in} \ (0,\infty) \times \Omega.
\end{equation}
Observe that in~\eqref{Ch1laecnolo} it appears the nonlocal sum of Laplacians $\Delta_{m^1\ast m^2}-\Delta_{m^1}-\Delta_{m^2}$ playing the role of the (local) fourth-order operator   $\Delta^2$ that appears in~\eqref{Ch1laec}.

 Formulation~\eqref{1901} allows us to study the (nonlocal) Cahn-Hilliard equation in a quite different   way than those used for the local and nonlocal  one, concretely we study it as  a  Lipschitz perturbation of the generalized porous medium equation, $$u_t = \Delta_{m^1} v,\
u  \in \gamma(v),$$ and we can do it for a large class of potentials that include the mentioned   ones. At our knowledge the approach we use  is new even for convolution kernels.
After the preliminaries we dedicate a section to   give   an overview of  the generalized nonlocal porous medium equation $$u_t = \Delta_{m^1} v,\
u  \in \gamma(v),$$ studied in~\cite{ST}, where, besides the nonlinearity driven by $\gamma$, general Leray-Lions type diffusion operators are considered. We will do it  for the sake  of completeness since it will serve as the basis to study the existence and uniqueness of mild and strong solutions for~\eqref{1901}.   We also   study the generalized nonlocal porous medium equation as a gradient flow in a \hbox{\lq\lq discrete $H^{-1}$-space"}, that is the nonlocal version of the results by Brezis \cite{Brezis0} for the local porous medium equation, which allows to get strong solutions in a  very direct way. These results are new in this abstract setting of random walk spaces.

It is worthily to mention that we cover the study of the Cahn-Hilliard system stated on a subset $\Omega\in \mathcal{B}$ with $0<\nu(\Omega)<\infty$ and $m_i$-connected, under (regional) homogeneous boundary conditions,
{\small \begin{equation}\label{DNLS1pre}
\left\{\begin{array}{lll} \displaystyle u_t(t,x) = \int_\Omega  (\mu(t,y)-\mu(t,x))d(m^1)_x(y) ,   &(t,x) \in (0, \infty) \times \Omega, \\[14pt]
\displaystyle \mu(t,x)= - \int_\Omega  (u(t,y)-u(t,x))d(m^2)_x(y)  + v(t,x)-cu(x,t),   &(t,x) \in (0, \infty) \times \Omega,
\\[14pt]
u(t,x)\in \gamma(v(t,x)), &(t,x) \in (0, \infty) \times \Omega,
\\[10pt] u(0,x) = u_0(x),   &x \in  \Omega,\end{array} \right.
\end{equation}}
\\ where  $\gamma$ is a maximal monotone graph with $0\in\gamma(0)$ and $c>0$, and, for which the following problem is a particular case:
{\small \begin{equation}\label{DNLS1preex}
\left\{\begin{array}{lll} \displaystyle u_t(t,x) = \int_\Omega  J_1(x-y)(\mu(t,y)-\mu(t,x))dy ,   &(t,x) \in (0, \infty) \times \Omega, \\[14pt]
\displaystyle \mu(t,x)= - \int_\Omega  J_2(x-y)(u(t,y)-u(t,x))dy  + v(t,x)-cu(x,t),   &(t,x) \in (0, \infty) \times \Omega,
\\[14pt]
u(t,x)\in \gamma(v(t,x)), &(t,x) \in (0, \infty) \times \Omega,
\\[10pt] u(0,x) = u_0(x),   &x \in  \Omega,\end{array} \right.
\end{equation}}
\\ where $\Omega$ is a bounded domain of $\mathbb{R}^N$ and $J_i:\mathbb{R}^N\to [0,+\infty[$, $i=1,2$, are non-singular  radial kernels with mass equal to 1. Taking into account, besides the kernels,  the general potentials considered  here, our existence and uniqueness results can be seen as a generalization of the  existence results for the weak-to-weak interaction   described by Gal and Shomberg in~\cite{GS}.

  We can put
  a    parameter $\delta>0 $ in the model (even one can deal  with $\delta=0$),
  \begin{equation}\label{DNLS1preforinipenalty}
\left\{\begin{array}{lll} \displaystyle u_t(t,x) = \Delta_{m^1}\mu(t,x) ,   &(t,x) \in (0, \infty) \times X, \\[10pt]
\displaystyle \mu(t,x)\in - \delta\Delta_{m^2}u(t,x) + \partial F(u),   &(t,x) \in (0, \infty) \times X,
\\[10pt] u(0,x) = u_0(x),   &x \in  X,\end{array} \right.
\end{equation}
   but for the existence results this is not relevant and we will assume $\delta=1$. We will consider such parameter $\delta$ in the last section where we study the asymptotic behaviour of the solutions. Then we will see that if (remember that we are dealing with $F$ like in~\eqref{dondelac})
   $$c<\delta \,\hbox{gap}(-\Delta_{m^2}),$$
   we have that for initial data $u_0 \in L^2$, the strong solution of \eqref{DNLS1preforinipenalty} converges, as time goes to infinity, to the media of $u_0$. This shows that, for $c$ small, this problem is not suitable for phase separation.   On the other hand, we will see that  for the   double-well obstacle potential, with $c \geq 2\delta$, the set of equilibria  contains the functions that separate the phases.

 Let us shortly describe the contents of the paper. In Section \ref{Pre} we recall all the notions about random walk spaces  required in this paper. Section \ref{PorousMedia} deals with the generalized porous medium equation in random walk spaces, that is one of the main tools used in the next section. We recall the results obtained in \cite{ST}, and prove some results by means of the Hilbertian theory that are   interesting by itself. In section \ref{Main} we obtain the main results about existence and uniqueness of solutions of the Cahn-Hilliard Equation.  First we study the Cauchy problem in $L^1$. As a consequence, we also obtain the existence and uniqueness of the regional Neumann problem. Moreover we   prove that the Cahn-Hilliard equation is the gradient flow of the Ginzburg-Landau free energy functional  in an adequate Hilbert space when the random walks have the same invariant measure. Finally, in Section~\ref{Properties} we obtain some properties of the solutions and  of their asymptotic behaviour, that particularly apply to   finite weighted discrete graphs and   non-smooth potentials.

\section{Preliminaries}\label{Pre}

 \subsection{Convex functions and subdifferentials}

 Let $H$ be a real Hilbert space with scalar product $\langle \cdot, \cdot \rangle_H$ and nor $\Vert u \Vert_H = \sqrt{\langle u, u \rangle_H}$. Given a function $\mathcal{F} : H \to ]-\infty, \infty]$ , we call the set $D(\mathcal{F}) : = \{ u \in H \ : \ \mathcal{F}(u) < + \infty \}$ the {\it effective domain} of $\mathcal{F}$, and $\mathcal{F}$ is said to be proper if $D(\mathcal{F})$ is non-empty.
Further, we say that $D(\mathcal{F})$ is {\it lower semi-continuous} if for every $c \in \R$, the sublevel set
$$E_c : = \{E c : = \{ u \in D(\mathcal{F}) \ : \ \mathcal{F}(u) \leq c \}$$
is closed in $H$.

 Given a   convex  proper function $\mathcal{F} : H \to ]-\infty, \infty]$, its {\it subdifferential}
 is defined by
 \begin{equation}\label{subdif2}\partial_H \mathcal{F} := \left\{(u,h) \in H \times H \ : \ \mathcal{F}(u+v) - \mathcal{F}(u) \geq \langle h, v \rangle_H \ \ \forall \, v \in D(\mathcal{F}) \right\}.
 \end{equation}
 This concept can be extended to any proper function $\mathcal{F} : H \to ]-\infty, \infty]$  by defining the {\it Gateaux subdifferential} of $\mathcal{F}$ as
 \begin{equation}\label{subdif1}\partial_H \mathcal{F} := \left\{ (u,h) \in H \times H \ : \ \liminf_{t \to 0^+} \frac{\mathcal{F} (u + tv) - \mathcal{F} (u)}{t} \geq \langle h, v \rangle_H \ \ \forall \, v \in D(\mathcal{F}) \right\},
 \end{equation}
 which,
 if $\mathcal{F}$ is convex, is reduced to the usual subdifferential.

  We say that $\alpha : [a,b] \rightarrow H$ is an {\it absolutely continuous curve}, and we write $\alpha \in AC([a,b]; H)$, if there exists $g \in L^1([a,b])$ such that
$$ \Vert \alpha(y) - \alpha(x) \Vert_H \leq \int_x^y g(t) dt \quad \forall \, a \leq x \leq y \leq b.$$

We define
$$AC_{loc}(0,\infty); H) := \{ \alpha : [0, \infty[ \rightarrow H \ : \ \alpha \in AC([a,b]; H) \ \hbox{for all $a < b$} \}.$$

 The following result is given in \cite[Proposition 11.4]{ABS}.

\begin{proposition} Let $I$ be an open interval of $\R$. Then any $u \in AC(I;H)$ is differentiable $\mathcal{L}^1$-a.e. $t \in I$, $\Vert u' \Vert \in L^1(I)$ and the fundamental theorem of calculus
$$u(t) - u(s) = \int_s^t u'(r) dr \quad \forall \, s,t \in I$$
holds.
\end{proposition}

The following definitions are given in \cite{ABS}

\begin{definition}\label{GF}{\rm Given a proper function $\mathcal{F} : H \to ]-\infty, \infty]$, we say that $u : (0, \infty) \rightarrow D(\mathcal{F} )$ is a  {\it a gradient flow} of  $\mathcal{F}$ if $u \in AC_{loc}(0,\infty); H)$ and
$$u'(t) +  \partial_H \mathcal{F} (u(t)) \ni 0 \quad \hbox{for $\mathcal{L}^1$-a.e. \ $t \in (0, \infty)$ }.$$
We say that $u$ starts from $u_0 \in H$ if $ \lim_{t \to 0} u(t) = u_0.$

}
\end{definition}
%
%
%

We are now in position to state the celebrated {\it Brezis-Komura Theorem} (see \cite{Brezis}, \cite{Komura}, or~\cite{ABS}).

\begin{theorem}\label{BKTheorem} Let $\mathcal{F} : H \to ]-\infty, \infty]$ be a proper convex and lower semi-continuous functional. Given $f \in L^2(0, T; H)$ and $u_0 \in \overline{D(\mathcal{F})}$ there exists a unique strong solution of the abstract Cauchy problem
\begin{equation}\label{ACHP}u'(t) +  \partial_H \mathcal{F} (u(t)) \ni f(t) \quad \hbox{for $\mathcal{L}^1$-a.e. \ $t \in (0, \infty)$ }\end{equation}
such that $u(0) =u_0$, that is $u \in C([0, + \infty[; H) \cap AC_{loc}((0,\infty); H)$ and satifies \eqref{ACHP}.

In the case $f =0$, if we denote $S(t) u_0 := u(t)$, the unique strong solution of the abstract Cauchy problem \eqref{ACHP}, then $S(t) : \overline{D(\mathcal{F})} \rightarrow H$ is a continuous semigroup satisfying the $T$-contraction property
$$\Vert (S(t) u_0 - S(t) v_0)^{\pm} \Vert \leq \Vert (u_0 -  v_0)^{\pm} \Vert \quad \forall \, t >0, \ \ u_0, v_0 \in \overline{D(\mathcal{F})}.$$
\end{theorem}

 \subsection{Random walk spaces}\label{RWS1} We recall some concepts and results about random walk spaces  given in \cite{MST0},  \cite{MST2} and \cite{MSTBook}.

 Let $(X,\mathcal{B})$ be a measurable space such that the $\sigma$-field $\mathcal{B}$ is countably generated.
A random walk $m$
on $(X,\mathcal{B})$ is a family of probability measures $(m_x)_{x\in X}$
on $\mathcal{B}$ such that $x\mapsto m_x(B)$ is a measurable function on $X$ for each fixed $B\in\mathcal{B}$.

The notation and terminology chosen in this definition comes from Ollivier's paper \cite{O}. As noted in that paper, geometers may think of $m_x$ as a replacement for the notion of balls around $x$, while in probabilistic terms we can rather think of these probability measures as defining a Markov chain whose transition probability from $x$ to $y$ in $n$ steps is
\begin{equation}
\displaystyle
dm_x^{*n}(y):= \int_{z \in X}  dm_z(y)dm_x^{*(n-1)}(z), \ \ n\ge 1
\end{equation}
and $m_x^{*0} = \delta_x$, the dirac measure at $x$.

\begin{definition}\label{convolutionofameasure}{\rm
If $m$ is a random walk on $(X,\mathcal{B})$ and $\mu$ is a $\sigma$-finite measure on $X$. The convolution of $\mu$ with $m$ on $X$ is the measure defined as follows:
$$\mu \ast m (A) := \int_X m_x(A)d\mu(x)\ \ \forall A\in\mathcal{B},$$
which is the image of $\mu$ by the random walk $m$.}
\end{definition}

\begin{definition}\label{def.invariant.measure} {\rm If $m$ is a random walk on $(X,\mathcal{B})$,
a $\sigma$-finite measure $\nu$ on $X$ is {\it invariant}
with respect to the random walk $m$ if
 $$\nu\ast m = \nu.$$

The measure $\nu$ is said to be {\it reversible} if moreover, the detailed balance condition $$dm_x(y)d\nu(x)  = dm_y(x)d\nu(y) $$ holds true.}
\end{definition}

\begin{definition}\label{DefMRWSf}{\rm
Let $(X,\mathcal{B})$ be a measurable space where the $\sigma$-field $\mathcal{B}$ is countably generated. Let $m$ be a random walk on $(X,\mathcal{B})$ and $\nu$ an invariant measure with respect to $m$. The measurable space together with $m$ and $\nu$ is then called a random walk space
and is denoted by $[X,\mathcal{B},m,\nu]$.}
\end{definition}

\begin{definition}\label{DefMRWS}{\rm
Let $[X,\mathcal{B},m,\nu]$ be a random walk space. If $(X,d)$ is a Polish metric space (separable completely metrizable topological space), $\mathcal{B}$ is its Borel $\sigma$-algebra and $\nu$ is a Radon measure (i.e. $\nu$ is inner regular
and locally finite) and we denote it by $[X,d,m,\nu]$.
}
\end{definition}

\begin{definition}\label{def.m.connected.random.walk.space}{\rm
Let $[X,\mathcal{B},m,\nu]$ be a random walk space. We say that $[X,\mathcal{B},m,\nu]$ is $m$-connected 
\index{maaconnected@$m$-connected!space} 
if, for every $D\in \mathcal{B}$ with $\nu(D)>0$ and $\nu$-a.e. $x\in X$,
$$\sum_{n=1}^{\infty}m_x^{\ast n}(D)>0.$$
}
\end{definition}

 \begin{definition}\label{def.m.interaction}{\rm
Let  $[X,\mathcal{B},m,\nu]$ be a random walk space and let $A$, $B\in\mathcal{B}$. We define the {\it $m$-interaction} between $A$ and $B$ as
\begin{equation}\label{m.interaction} L_m(A,B):= \int_A \int_B dm_x(y) d\nu(x)=\int_A m_x(B) d\nu(x).
 \end{equation}
 }
 \end{definition}

  The following result gives a characterization of $m$-connectedness in terms of the $m$-interaction between sets.

\begin{proposition}\label{connectedness.iff.Lm}(\cite[Proposition 2.11]{MST0}, \cite[Proposition 1.34]{MSTBook})
 Let $[X,\mathcal{B},m,\nu]$ be a random walk space.
The following statements are equivalent:
\item{ (i) } $[X,\mathcal{B},m,\nu]$ is $m$-connected.
\item {(ii)} If $ A,B\in\mathcal{B}$ satisfy $A\cup B=X$ and $L_m(A,B)= 0$, then either $\nu(A)=0$ or $\nu(B)=0$.
\item {(iii)} If $A\in \mathcal{B}$ is a $\nu$-invariant set then either $\nu(A)=0$ or $\nu(X\setminus A)=0$.
\end{proposition}

Let us see now some examples of random walk spaces.

 \begin{example}\label{example.nonlocalJ} \rm
Consider the metric measure space $(\R^N, d, \mathcal{L}^N)$, where $d$ is the Euclidean distance and $\mathcal{L}^N$ the Lebesgue measure on $\R^N$ (which we will also denote by $|.|$). For simplicity, we will write $dx$ instead of $d\mathcal{L}^N(x)$. Let  $J:\R^N\to[0,+\infty[$ be a measurable, nonnegative and radially symmetric
function  verifying  $\int_{\R^N}J(x)dx=1$. Let $m^J$ be the following random walk on $(\R^N,d)$: 
\index{m-J@$m^J$} 
$$m^J_x(A) :=  \int_A J(x - y) dy \quad \hbox{ for every $x\in \R^N$ and every Borel set } A \subset  \R^N  .$$
Then, applying Fubini's Theorem it is easy to see that the Lebesgue measure $\mathcal{L}^N$ is reversible with respect to $m^J$. Therefore, $[\R^N, d, m^J, \mathcal{L}^N]$ is a reversible metric random walk space.
\end{example}

\begin{example}\label{example.graphs}[Weighted discrete graphs] \rm Consider a locally finite  weighted discrete graph $$G = (V(G), E(G)),$$ where $V(G)$ is the vertex set, $E(G)$ is the edge set and each edge $(x,y) \in E(G)$ (we will write $x\sim y$ if $(x,y) \in E(G)$) has a positive weight $w_{xy} = w_{yx}$ assigned. Suppose further that $w_{xy} = 0$ if $(x,y) \not\in E(G)$.  Note that there may be loops in the graph, that is, we may have $(x,x)\in E(G)$ for some $x\in V(G)$ and, therefore, $w_{xx}>0$. Recall that a graph is locally finite if every vertex is only contained in a finite number of edges.

 A finite sequence $\{ x_k \}_{k=0}^n$  of vertices of the graph is called a {\it  path} if $x_k \sim x_{k+1}$ for all $k = 0, 1, ..., n-1$. The {\it length} of a path $\{ x_k \}_{k=0}^n$ is defined as the number $n$ of edges in the path. With this terminology, $G = (V(G), E(G))$ is said to be {\it connected} if, for any two vertices $x, y \in V$, there is a path connecting $x$ and $y$, that is, a path $\{ x_k \}_{k=0}^n$ such that $x_0 = x$ and $x_n = y$.  Finally, if $G = (V(G), E(G))$ is connected, the {\it graph distance} $d_G(x,y)$ between any two distinct vertices $x, y$ is defined as the minimum of the lengths of the paths connecting $x$ and $y$. Note that this metric is independent of the weights.

For $x \in V(G)$ we define the weight at $x$ as
$$d_x:= \sum_{y\sim x} w_{xy} = \sum_{y\in V(G)} w_{xy},$$
and the neighbourhood of $x$ as $N_G(x) := \{ y \in V(G) \, : \, x\sim y\}$. Note that, by definition of locally finite graph, the sets $N_G(x)$ are finite. When all the weights are $1$, $d_x$ coincides with the degree of the vertex $x$ in a graph, that is,  the number of edges containing $x$.

For each $x \in V(G)$  we define the following probability measure 
  \index{m-G@$m^G$} 
\begin{equation}\label{discRW}m^G_x:=  \frac{1}{d_x}\sum_{y \sim x} w_{xy}\,\delta_y.\\ \\
\end{equation}
It is not difficult to see that the measure $\nu_G$ defined as
 $$\nu_G(A):= \sum_{x \in A} d_x,  \quad A \subset V(G),$$
is a reversible measure with respect to this random walk. Therefore, $[V(G),\mathcal{B},m^G,\nu_G]$ is a reversible random walk space ($\mathcal{B}$ is the $\sigma$-algebra of all subsets of $V(G)$) and $[V(G),d_G,m^G,\nu_G]$ is a reversible metric random walk space.

In Machine Learning Theory (\cite{G-TS}, \cite{G-TSBLBr}), an example of a weighted discrete graph is a point cloud in $\R^N$,  $V=\{x_1, \ldots , x_n \}$,  with edge weights $w_{x_i,x_j}$ given by
$$w_{x_i,x_j} := \eta(\vert x_i - x_j\vert), \quad 1 \leq i,j \leq n,$$
where the kernel  $ \eta : [0, \infty) \rightarrow [0, \infty)$ is a radial profile satisfying
\begin{itemize}
\item[(i)] \ $ \eta(0) >0$, and $\eta$ is continuous at $0$, \vspace{2pt}
\item[(ii)]  \ $\eta$ is non-decreasing, \vspace{2pt}
\item[(iii)] \ and the integral $\int_0^\infty  \eta(r) r^N dr$ is finite.
\end{itemize}
\end{example}

\begin{example}\label{example.markov.kernel} [Markov chains] \rm Let $K: X \times X \rightarrow \R$ be a Markov kernel on a countable space $X$, i.e.,
 $$K(x,y) \geq 0 \quad \forall x,y \in X, \quad \quad \sum_{y\in X} K(x,y) = 1 \quad \forall x \in X.$$

Then, if 
\index{m-K@$m^K$} 
 $$m^K_x(A):= \sum_{y \in A} K(x,y), \ \ x\in X, \, A\subset X$$
 and $\mathcal{B}$ is the $\sigma$-algebra of all subsets of $X$, $m^K$ is a random walk on $(X,\mathcal{B})$.

  Recall that, in discrete Markov chain theory terminology, a measure $\pi $ on $X$ satisfying
 $$\sum_{x \in X} \pi(x) = 1 \quad \hbox{and} \quad \pi(y) = \sum_{x \in X} \pi(x) K(x,y)  \quad \forall y \in X,$$
 is called a stationary probability measure (or steady state) on $X$. Of course, $\pi$ is a stationary probability measure if, and only if, $\pi$ is and invariant probability measure with respect to $m^K$. Consequently, if $\pi$ is a stationary probability measure on $X$, then $[X, \mathcal{B}, m^K, \pi]$ is a random walk space.

Furthermore, a stationary probability measure $\pi$ is said to be reversible for $K$ if the following detailed balance equation holds:
 $$K(x,y) \pi(x) = K(y,x) \pi(y) \ \hbox{ for } x, y \in X.$$
 This balance condition is equivalent to
 $$dm^K_x(y)d\pi(x)  =   dm^K_y(x)d\pi(y)\ \hbox{ for } x, y \in X.$$

Note that, given a locally finite weighted discrete graph $G = (V(G), E(G))$ as in Example \ref{example.graphs}, there is a natural definition of a Markov chain
on the vertices. Indeed, define the Markov kernel  $K_G: V(G)\times V(G) \rightarrow \R$ as
$$K_G(x,y):= \frac{1}{d_x}  w_{xy}.$$
Then, $m^G$ and $m^{K_G}$ define the same random walk.  If $\nu_G(V(G))$ is finite, the unique reversible probability measure with respect to $m^G$ is given by
$$\pi_G(x):= \frac{1}{\nu_G(V(G))} \sum_{z \in V(G)} w_{xz}. $$
\end{example}

\begin{example}\label{example.restriction.to.Omega} \rm Given a random walk  space $[X,\mathcal{B},m,\nu]$ and $\Omega \in \mathcal{B}$ with $\nu(\Omega) > 0$, let 
\index{m-Omega@$m^\Omega$} 
$$m^{\Omega}_x(A):=\int_A d m_x(y)+\left(\int_{X\setminus \Omega}d m_x(y)\right)\delta_x(A) \quad \hbox{ for every } A\in\mathcal{B}_\Omega  \hbox{ and } x\in\Omega.
$$
Then, $m^{\Omega}$ is a random walk on $(\Omega,\mathcal{B}_\Omega)$ and it easy to see that $\nu \res \Omega$ is invariant with respect to $m^{\Omega}$. Therefore,  $[\Omega,\mathcal{B}_\Omega,m^{\Omega},\nu \res \Omega]$ is a random walk space. Moreover, if $\nu$ is reversible with respect to $m$ then $\nu \res \Omega$ is  reversible with respect to $m^{\Omega}$. Of course, if $\nu$ is a probability measure we may normalize $\nu \res \Omega$ to obtain the random walk space
$$\left[\Omega,\mathcal{B}_\Omega,m^{\Omega}, \frac{1}{\nu(\Omega)}\nu \res \Omega \right].$$
 Note that, if $[X,d,m,\nu]$ is a metric random walk space and $\Omega$ is closed, then $[\Omega,d,m^{\Omega},\nu \res \Omega]$ is also a metric random walk space, where we abuse notation and denote by $d$ the restriction of $d$ to $\Omega$.

In particular, in the context of Example \ref{example.nonlocalJ}, if $\Omega$ is a closed and bounded subset of $\R^N$, we obtain the metric random walk space $[\Omega, d, m^{J,\Omega},\mathcal{L}^N\res \Omega]$ where 
\index{m-JOmega@$m^{J,\Omega}$} 
$m^{J,\Omega} := (m^J)^{\Omega}$; that is,
$$m^{J,\Omega}_x(A):=\int_A J(x-y)dy+\left(\int_{\R^n\setminus \Omega}J(x-z)dz\right)d\delta_x$$ for every Borel set   $A \subset  \Omega$  and $x\in\Omega$.

\end{example}

\subsection{The nonlocal gradient, divergence and Laplace operators}\label{nonlocal.notions.1.section}

Let us introduce the nonlocal counterparts of some classical concepts.

\begin{definition}\label{nonlocalgraddiv}{\rm
Let $[X,\mathcal{B},m,\nu]$ be a random walk space. Given a function $u : X \rightarrow \R$ we define its {\it nonlocal gradient}
$\nabla u: X \times X \rightarrow \R$ as
$$\nabla u (x,y):= u(y) - u(x) \quad \forall \, x,y \in X.$$
Moreover, given $\z : X \times X \rightarrow \R$, its {\it $m$-divergence} 
${\rm div}_m \z : X \rightarrow \R$ is defined as
 $$({\rm div}_m \z)(x):= \frac12 \int_{X} (\z(x,y) - \z(y,x)) dm_x(y).$$
 }
\end{definition}

\begin{definition}\label{def.averaging.operator}{\rm
If $\nu$ is an invariant measure with respect to $m$, we define the linear operator $M_m$ on $L^1(X,\nu)$ into itself as follows
$$M_m f(x):=\int_X f(y)dm_x(y), \ \ f\in L^1(X,\nu).$$
$M_m$ is called the {\it averaging operator} 
on $[X,\mathcal{B},m]$}
\end{definition}
Note that, if $f\in L^1(X,\nu)$ then, using the invariance of $\nu$ with respect to $m$,
$$\int_X\int_X |f(y)|dm_x(y)d\nu(x)=\int_X |f(x)| d\nu(x) < \infty,$$
so $f\in L^1(X,m_x)$ for $\nu$-a.e. $x\in X$, thus $M_m$ is well defined from $L^1(X,\nu)$ into itself.

\begin{remark}\label{remark.Mm} \rm Let $\nu$ be an invariant measure with respect to $m$. It follows that
$$\Vert M_m f\Vert_{L^1(X,\nu)}\leq \Vert f\Vert_{L^1(X,\nu)} \ \ \forall f\in L^1(X,\nu),$$
so that $M_m$ is a contraction on $L^1(X,\nu)$. In fact, since $M_m f\ge 0$ if $f\ge 0$, we have that $M_m$ is a positive contraction on $L^1(X,\nu)$.

Moreover, by Jensen's inequality, we have that, $$ \Vert M_m f \Vert^2_{L^2(X,\nu)}\leq  \int_X f^2(x) d\nu(x) =\Vert f \Vert^2_{L^2(X,\nu)}.
$$
Therefore, $M_m$ is a linear operator in $L^2(X,\nu)$  with  domain
$$D(M_m)= L^1(X, \nu) \cap  L^2(X, \nu).$$
Consequently, if $\nu(X)<+\infty$, $M_m$ is a bounded linear operator from $L^2(X,\nu)$ into itself satisfying $\Vert M_m  \Vert=\Vert M_m  \Vert_{\mathcal{B}(L^2(X,\nu),L^2(X,\nu))}  \leq 1$.  $\blacksquare$
\end{remark}

We define the (nonlocal) Laplace operator as follows.
\begin{definition}\label{deflap1310}{\rm
Let $[X,\mathcal{B},m,\nu]$ be a random walk space, we define the {\it $m$-Laplace operator} (or {\it $m$-Laplacian}) from $L^1(X,\nu)$ into itself as $\Delta_m:= M_m - I$, i.e.,
$$\Delta_m u(x)= \int_X u(y) dm_x(y) - u(x) = \int_X (u(y) - u(x)) dm_x(y), \ \ u\in L^1(X,\nu).$$}
\end{definition}
 Note that
$$\Delta_m f (x) = {\rm div}_m (\nabla f)(x).$$

\begin{remark}\label{remark.laplacian} \rm
We have that  $\Vert \Delta_mf\Vert_1\le 2\Vert f\Vert_1$ and
\begin{equation}\label{Lap0}
\int_X \Delta_m f(x) d\nu(x) = 0 \quad \hbox{$\forall\, f\in L^1(X,\nu)$}.
\end{equation}
As in Remark \ref{remark.Mm} we obtain that $\Delta_m$ is a linear operator in $L^2(X,\nu)$  with  domain $$D(\Delta_m) = L^1(X, \nu) \cap  L^2(X, \nu).$$
Moreover,   if $\nu(X)<+\infty$,  $\Delta_m$ is a bounded  linear operator in $L^2(X,\nu)$  satisfying $\Vert \Delta_m \Vert  \leq 2$.  $\blacksquare$
\end{remark}

We define the energy functional $\mathcal{H}_m : L^2(X, \nu) \rightarrow [0, + \infty]$  defined as
$$\mathcal{H}_m(f)= \left\{ \begin{array}{ll} \displaystyle\frac{1}{4} \int_{X \times X} (f(x) - f(y))^2 dm_x(y) d\nu(x) \quad &\hbox{ if $f\in L^2(X, \nu) \cap  L^1(X, \nu)$.} \\[12pt] + \infty, \quad &\hbox{else}. \end{array}\right.$$
We denote
$$D(\mathcal{H}_m)=L^2(X, \nu) \cap  L^1(X, \nu).$$
In \cite{MST0} it is proved that
$$-\Delta_m = \partial_{L^2(X, \nu)} \mathcal{H}_m.$$

In the case of the random walk space associated with a locally finite weighted discrete graph $G=(V,E)$ (as defined in Example~\ref{example.graphs}), the $m^G$-Laplace operator coincides with the graph Laplacian (also called the normalized graph Laplacian) studied by many authors (see, for example, \cite{BJ}, \cite{BJL}, \cite{DK}, \cite{Elmoatazetal}, \cite{Hafiene} or  \cite{JL}):
$$\Delta u(x):=\frac{1}{d_x}\sum_{y\sim x}w_{xy}(u(y)-u(x)), \quad u\in L^2(V,\nu_G), \ x\in V .$$

\begin{proposition}\label{integration.by.parts}(Integration by parts formula)
Let $[X,\mathcal{B},m,\nu]$ be a reversible random walk space. Then, 
\index{integration by parts formula} 
 \begin{equation}\label{intbpart}
 \int_X f(x) \Delta_m g (x) d\nu(x) =  -\frac {1}{2} \int_{X \times X} \nabla f(x,y)\nabla g(x,y) d(\nu\otimes m_x)(x,y)
 \end{equation}
 for $f,g \in L^1(X, \nu)\cap  L^2(X, \nu)$. In particular
  for $f\in D(\mathcal{H}_m)$, we have
$$ \mathcal{H}_m(f) =   - \frac12\int_X f(x) \Delta_m f (x) d\nu(x).
$$
\end{proposition}

\begin{definition}\label{defpoin}{\rm
We say that  $[X,\mathcal{B},m,\nu]$  satisfies a {\it Poincar\'{e} inequality} if there exists $\lambda >0$ such that
\begin{equation}\label{Poinca1}\lambda\Vert f\Vert_{L^2(X, \nu)}^2\le \mathcal{H}_m(f)\quad  \hbox{for all} \ f \in L^2(X,\nu) \hbox{ with } \int_X f d\nu=0.\end{equation}
}
\end{definition}

In \cite{MST0} (see also \cite{MSTBook}) it is shown that under quite general assumptions such an inequality holds true for the examples of random walk spaces given in Subsection \ref{RWS1}.

  The {\it spectral gap} of $-\Delta_m$ is defined as
\begin{equation}\label{ladefdegap}{\rm gap}(-\Delta_m) := \inf \left\{ \frac{2\mathcal{H}_m(f)}{\Vert f \Vert^2_2 } \ : \ f \in D(\mathcal{H}_m), \ \Vert f \Vert_2 \not= 0, \ \int_X f d\nu = 0 \right\}.
\end{equation}
We have that $\frac12 {\rm gap}(-\Delta_m)$ is the best constant in the Poincar\'{e} inequality.   Moreover, it is well-known that ${\rm gap}(-\Delta_m)\le 2$.

In \cite{MST2} (see also  \cite{MSTBook})  we introduce and study the following concept that will be used later on.

\begin{definition}\label{Mcurv}\rm  Let $[X,\mathcal{B},m,\nu]$ be a random walk space and let $E \subset X$ be $\nu$-measurable. For a point $x  \in X$ we define its {\it $m$-mean curvature} as
\begin{equation}\label{defcur}\mathcal{H}^m_{\partial E}(x):= \int_{X}  (\1_{X \setminus E}(y) - \1_E(y)) dm_x(y) = 1 - 2 m_x(E),\end{equation}
 which takes values in $[-1,1]$.
\end{definition}

\section{The Cauchy problem for the generalized porous medium equation}\label{PorousMedia}

  Assume that $[X,\mathcal{B},m,\nu]$ is a  $m$-connected random walk space  with $\nu$ reversible, and   $\nu(X) < +\infty$. In this section we will study the Cauchy problem for the generalized porous medium equation under two points of view.

Let $\gamma$ be a maximal monotone graph in $\R \times \R$ such that $0 \in \gamma(0)$, and set
  $$ \gamma^-:=\inf\{\mbox{Ran}(\gamma)\},\ \gamma^+:=\sup\{\mbox{Ran}(\gamma)\}.$$
 We assume $\gamma^-<\gamma^+$.    Let us also define   $$j_\gamma^*(r):=\displaystyle\int_0^r(\gamma^{-1})^0(s)ds,$$  for $r\in D(\gamma^{-1})$ (the effective domain of $\gamma^{-1}$),  where $(\gamma^{-1})^0(r)$ is the element of $\gamma^{-1}(r)$ of least norm. We have that
 $$\partial j_\gamma^*(r)=\gamma^{-1}(r).$$

 Consider the problem
\begin{equation}\label{PME2}
 \left\{ \begin{array}{ll}\displaystyle
 \frac{\partial u}{\partial t} - \Delta_m v = f \quad &\hbox{in} \ \ (0, T) \times X, \\[12pt]
 u\in\gamma(v)\quad &\hbox{in} \ \ (0, T) \times X,
 \\[10pt]
 u(0, x) = u_0(x) \quad & x \in X.\end{array} \right.
 \end{equation}
  Problem \eqref{PME2}, for different choices of $\gamma$, gives rise to important examples. For instance, if  $\gamma^{-1}(r) = \vert r \vert^{m-1} r$, and $m >1$, it corresponds to the (nonlocal) {\it  porous medium equation}, if  $m=1$ is the (nonlocal)
 heat equation, and if $0 < m < 1$, it is the {\it fast diffusion equation}. For
$$\gamma^{-1}(r)= \left\{ \begin{array}{lll} r \quad &\hbox{if} \ \ r < 0, \\[4pt][0,1] \quad &\hbox{if} \ \ 0\leq r \leq 1, \\[4pt] r - 1 \quad &\hbox{if} \ \ r \geq 1, \end{array} \right. $$
it deals with a (nonlocal) {\it Stephan problem}.  While (nonlocal) Hele-Shaw type problems correspond to a choice like
  $$\gamma(r):=  \left\{\begin{array}{ll}
0 & \hbox{if $r<0$},\\[4pt]
[0,1]&\hbox{if $r=0$,}\\[4pt]
1  & \hbox{if $r>0$}.
\end{array}\right.$$

 As we will see another important examples of $\gamma$ are the related with the graphs  that appear in the Ginzburg-Landau free energy  for the Cahn-Hilliard system.

\subsection{The $L^1$-Theory}

  Let $T>0$ and  $f\in L^1(0,T;L^1(X,\nu))$.
 The study of  existence and uniqueness of solution of Problem \eqref{PME2} was done in~\cite{ST} by means of the Nonlinear Semigroup Theory  (\cite{Cr}, \cite{Cr2}, \cite{CrandallLiggett}, \cite{Barbu}).  Problem \eqref{PME2} is written as an abstract Cauchy problem in  $L^1(X, \nu)$ associated with a T-accretive operator.

 \begin{definition}\label{operator1}{\rm Define in $L^1(X, \nu)$ the operator ${\mathbf B}^m_\gamma$ as
$(u, \hat{u}) \in {\mathbf B}^m_\gamma$  if: \\[6pt]
$u,\hat{u} \in L^1(X, \nu)$ and there exists $v \in L^2(X, \nu)$ with
$$u \in \gamma(v) \ \hbox{$\nu$-a.e.}$$
such that
$$- \Delta_{m} v=\hat{u}.$$
}
 \end{definition}

With this operator at hand, Problem~\eqref{PME2} can be rewritten as the following abstract Cauchy problem:
\begin{equation}\label{ACProbPM}
\left\{ \begin{array}{ll} u'(t) +  {\mathbf B}^m_\gamma  (u(t)) \ni f(t),\quad t>0, \\[8pt] u(0) = u_0.
\end{array} \right.
\end{equation}
 The following facts are proved in \cite[Theorems 3.2, 3.3, 3.4]{ST} (for, using the notation of such reference, $\gamma=\beta$, $a_p(x,y,r)=r$ and $\Omega_1\cup\Omega_2=X$):

 \begin{proposition}[\cite{ST}]\label{accretive}  Assume $[X,\mathcal{B},m,\nu]$  satisfies a {\it Poincar\'{e} inequality}.
  Then:\\
 1.The domain of the operator ${\mathbf B}^m_\gamma$ satisfies:
 $$\overline{D({\mathbf B}^m_\gamma)}^{L^1(\Omega)}=\left\{u\in L^1(X,\nu): \gamma^-\le u \le \gamma^+ \right\}.$$
 2. ${\mathbf B}^m_\gamma$  is  {\rm T}-accretive in $L^1(X, \nu)$ and satisfies the range condition:
\begin{equation}\label{e1RCond1}
\left\{ u \in L ^2(X, \nu)   :   \nu(X) \gamma^- < \int_X u d\nu  < \nu(X) \gamma^+ \right\}  \subset R(I + \lambda {\mathbf B}^m_\gamma) \quad\forall\lambda>0.
\end{equation}
3.  For any $T>0$, and for any $u_0\in \overline{D({\mathbf B}^m_\gamma)}^{L^1(\Omega)}$ and $f\in L^1(0,T;L^1(\Omega,\nu))$ satisfying   \begin{equation}\label{natcompcond}\displaystyle
\nu(X) \gamma^- < \int_X u_0 d\nu
+ \int_0^t\int_X f d\nu dt < \nu(X) \gamma^+, \quad\forall 0\le t\le T,
\end{equation}
 there exists a unique mild-solution $u\in C([0,T]:L^1(\Omega,\nu))$ of Problem~\eqref{ACProbPM} (hence of Problem~\eqref{PME2}).

4. Let $u_0,\widetilde{u_0}\in \overline{D({\mathbf B}^m_\gamma)}^{L^1(\Omega)}$ and $f,\widetilde{f}\in L^1(0,T;L^1(\Omega,\nu))$, satisfying the corresponding assumption~\eqref{natcompcond}, and  $u,\widetilde{u}$ the respective mild solutions of   Problem~\eqref{ACProbPM}, then
\begin{equation}\label{1553}
\begin{array}{ll}
\displaystyle \int_X (u(t,x)-\widetilde u(t,x))^+d\nu(x)\le \int_X(u_0(x)-\widetilde u_0(x))^+d\nu(x)+
\\ \\
\displaystyle \phantom{\int_X (v(t,x)-\widetilde u(t,x))^+d\nu(x)\le} +\int_0^t\int_X\left(f(s,x)-\widetilde f(s,x)\right)^+\nu(x)ds,\quad\forall 0\le t\le T.
\end{array}
\end{equation}

5. If in addition
$u_0\in L^2(\Omega,\nu)$ and $\displaystyle\int_Xj_\gamma^*(u_0)d\nu<+\infty$, and $f\in L^2(0,T; L^2(X,\nu))$,  the mild solution  belongs to $W^{1,1}(0,T;L^2(\Omega))$, therefore it is  a strong solution.
\end{proposition}

  In the next section we see that  strong solutions  can be also obtained via a gradient flow in an adequate Hilbert space. This will be used later on to get strong solutions for the Cahn-Hilliard problem.

\subsection{The Hilbertian theory} Now we study the nonlocal version of the results by Brezis \cite{Brezis0} for the local porous medium equation.
In order to do this let us introduce a Hilbertian structure in random walk spaces.   The particular case of finite weighted graphs and $\gamma^{-1}$ an increasing function was consider  in \cite{EM} where the discrete porous medium arise as gradient flow of certain entropy functionals with respect to suitable non-local transportation metrics.

Let $[X,\mathcal{B},m,\nu]$ be a $m$-connected random walk space with $\nu$ reversible,  $\nu(X) < +\infty$, and assume that $[X,\mathcal{B},m,\nu]$  satisfies a Poincar\'{e} inequality.

 We have that $\Delta_m$ is a linear bounded operator in $L^2(X, \nu)$, and by the ergodicity of $\nu$, we have
 $${\rm Ker} (\Delta_m) = {\rm Lin} \{ \1_X \}.$$
 Since $\Delta_m$ is selfadjoint in $L^2(X, \nu)$ and $[X,\mathcal{B},m,\nu]$ satisfies a  Poincar\'{e} inequality, we have that
 $${\rm Ran}(\Delta_m) = L_0^2(X, \nu):=\left \{ u \in L^2(X, \nu) \ : \ \int_X u(x) d \nu(x) = 0 \right\}$$
 and
 moreover $$\Delta_m:L_0^2(X, \nu)\to L_0^2(X, \nu)$$ is bijective. For $v\in  L_0^2(X, \nu)$, $\Delta_m^{-1} v$ denotes the preimage of $v$ via this bijection.

 From now on we will denote
 $$H_m^{-1}(X, \nu):= {\rm Ran}(\Delta_m)=L_0^2(X, \nu)$$
 endowed the inner product
 $$\langle v_1, v_2 \rangle_{H_m^{-1}}:= - \int_X \Delta_m^{-1} v_1 v_2 d\nu.$$
 By the integration by parts formula \eqref{intbpart}, we have
 $$\langle v_1, v_2 \rangle_{H_m^{-1}} = \frac12 \langle \nabla \phi_1, \nabla \phi_2 \rangle_{L^2(X \times X, d(\nu\otimes m_x))},$$
 being $\phi_1, \phi_2 \in H_m^{-1}$ the unique functions such that $\Delta_m \phi_i = v_i$, $i=1,2$., that is
 $$\langle v_1, v_2 \rangle_{H_m^{-1}} = \frac12 \langle \nabla \Delta_m^{-1} v_1, \nabla \Delta_m^{-1} v_2 \rangle_{L^2(X \times X, d(\nu\otimes m_x))}.$$  Let   $||.||_{H_m^{-1}}$ be the induced norm by such inner product.

\begin{proposition}\label{isomorph}  Assume  that the random walk space  $[X,\mathcal{B},m,\nu]$  satisfies a {\it Poincar\'{e} inequality}  (with constant $\lambda$).  Then,  the  Hilbert space $H_m^{-1}(X, \nu)$ is  isomorphic to the Hilbert space $(L_0^2(\Omega, \nu), \Vert \cdot \Vert_{L^2(\Omega, \nu)})$ with
$$\Vert v \Vert_{L^2(X, \nu)} \leq \sqrt{\frac{2}{\lambda}} \Vert v \Vert_{H^{-1}_m}  \leq  \frac{1}{\lambda} \Vert v \Vert_{L^2(X, \nu)}.$$
\end{proposition}
\begin{proof} Given $v \in H_m^{-1}(X, \nu)$, applying Poincar\'{e}'s inequality, we have
\begin{equation}\label{LL001}\Vert(\Delta_m^{-1} v) \Vert^2_{L^2(X, \nu)} \leq \frac{1}{\lambda} \mathcal{H}_m (\Delta_m^{-1} v)$$ $$ = \frac{1}{4\lambda}\int_{X \times X} \nabla \Delta_m^{-1} v \cdot \nabla \Delta_m^{-1} v d(\nu\otimes m_x)= \frac{1}{2\lambda} \Vert v \Vert^2_{H^{-1}_m}.\end{equation}
Then
$$\Vert v \Vert^2_{L^2(X, \nu)} = \Vert \Delta_m (\Delta_m^{-1} v) \Vert^2_{L^2(X, \nu)} \leq 4 \Vert \Delta_m^{-1} v \Vert^2_{L^2(X, \nu)} \leq\frac{2}{\lambda} \Vert v \Vert^2_{H^{-1}_m}.$$
Hence
$$\Vert v \Vert_{L^2(X, \nu)} \leq \sqrt{\frac{2}{\lambda}} \Vert v \Vert_{H^{-1}_m}.$$

On the other hand, applying  Cauchy-Scharz  and~\eqref{LL001}, we have
$$\Vert v \Vert^2_{H^{-1}_m} = - \int_X v \Delta_m^{-1} v d \nu \leq \Vert v \Vert_{L^2(X, \nu)} \Vert \Delta_m^{-1} v  \Vert_{L^2(X, \nu)} \leq \Vert v \Vert_{L^2(X, \nu)} \Vert \frac{1}{\sqrt{2\lambda}} \Vert v \Vert_{H^{-1}_m}.$$
Hence
\[\Vert v \Vert_{H^{-1}_m} \leq  \frac{1}{\sqrt{2\lambda}} \Vert v \Vert_{L^2(X, \nu)}.\qedhere\]
\end{proof}

 Our aim is to study  Problem~\eqref{PME2} as a gradient flow in the Hilbert space $H^{-1}_m(X, \nu)$. For this, we consider the energy functional
$\Psi_\gamma :  H^{-1}_m(X, \nu)\rightarrow ]\infty, +\infty]$ defined by
\begin{equation}\label{1744}\Psi_\gamma(u):= \left\{ \begin{array}{ll} \displaystyle\int_X j^*_\gamma(u(x)) d\nu(x) \quad &\hbox{if} \ \   j^*_\gamma(u) \in L^1(X, \nu), \\ \\ + \infty \quad &\hbox{otherwise}.\end{array} \right.
\end{equation}

\begin{theorem}\label{Charact1}
We have that $\Psi_\gamma$ is convex and lower semi-continuous on $H_m^{-1}(X, \nu)$. Moreover,
\begin{equation}\label{e2Chact}\begin{array}{ll}
\partial_{H^{-1}_m(X, \nu)} \Psi_\gamma=  &\Big\{ (u, w) \in H_m^{-1}(X, \nu) \times H_m^{-1}(X, \nu)   :   w = - \Delta_m v, \\[8pt] & \qquad v \in L^2(X, \nu), \ u(x)\in\gamma(v(x)) \ \nu-a.e. \ x \in X \Big\}.\end{array}
\end{equation}
 \end{theorem}

\begin{proof} Obviously  $\Psi_\gamma$ is convex.  On the other hand,  we have that $\Psi_\gamma$ is lower semi-continuous in $L^1(X, \nu)$ (see for example \cite[Proposition 2.7]{Barbu}), then by Proposition~\ref{isomorph},  we get that $\Psi_\gamma$ is lower semi-continuous on $H^{-1}(X, \nu)$.  Therefore, $\partial_{H^{-1}_m(X, \nu)} \Psi_\gamma$ is a maximal monotone operator.

We consider now the operator
$$\begin{array}{ll}
A_\gamma =  &\{ (u, w) \in H_m^{-1}(X, \nu) \times H_m^{-1}(X, \nu)   :   w = - \Delta_m v, \\ \\ & v \in L^2(X, \nu), \ u(x) \in \gamma (v(x)) \ \nu-a.e. \ x \in X \}.\end{array} $$
Let us see now that $\partial_{H^{-1}_m(X, \nu)} \Psi_\gamma=A_\gamma$. For that we will  see first that $A_\gamma \subset \partial_{H^{-1}_m(X, \nu)} \Psi_\gamma$ (this implies monotonicity of $A_\gamma$) and next that $A_\gamma$ is in fact maximal monotone.

Given  $(u, w) \in A_\gamma$, we have $w = - \Delta_m v$, with $v \in L^2(X, \nu), \ u(x) \in \gamma (v(x)) \ \nu-a.e. \ x \in X$. Then, $v(x) \in \gamma^{-1}(u(x))=\partial j^*_\gamma (u(x)) \ \nu-a.e. \ x \in X$, and consequently, if $\tilde{u} \in H^{-1}(X,\nu)$,
 $$j^*_\gamma(\tilde{u}(x)) - j^*_\gamma(u(x) \geq v(x)(\tilde{u}(x) - u(x)) \quad \ \nu-a,e. \ x \in X,$$ and hence
 $$\Psi_\gamma(\tilde{u}) - \Psi_\gamma(u) \geq \int_X v(x)(\tilde{u}(x) - u(x)) d \nu(x) = \int_X (- \Delta_m^{-1} w) (\tilde{u} - u) d \nu = \langle w, \tilde{u} - u \rangle_{H_m^{-1}}. $$
 Therefore $(u, w) \in \partial_{H_m^{-1}(X, \nu)} \Psi_\gamma$, and, consequently, $A_\gamma \subset \partial_{H_m^{-1}(X, \nu)} \Psi_\gamma$.

 By Minty's Theorem, to see that $A_\gamma$ is maximal  monotone, we need to show that verifies the range condition
\begin{equation}\label{e5Chact}
R(I + A_\gamma) = H_m^{-1}(X,\nu).
\end{equation}
 Then, we must show that given $z \in  H^{-1}(X,\nu)$, there exists $u \in H_m^{-1}(X,\nu)$ such that $(u, z- u) \in A_\gamma$.  Now, this is equivalent to show that there exists $u \in L^2_0(X,\nu)$ such that $(u, z - u) \in {\mathbf B}^m_\gamma\cap L^2_0(X,\nu)\times L^2_0(X,\nu)$ which is true, thanks to Proposition~\ref{accretive}.{\it 2.} Indeed, there exists $u \in L^1(X,\nu)$ such that $(u, z - u) \in {\mathbf B}^m_\gamma$, but it is easy that, in fact, $u\in L^2_0(X,\nu).$
\end{proof}

By the above theorem and the Brezis-Komura Theorem (Theorem \ref{BKTheorem}), we have:

\begin{theorem}\label{1612}
  For every inital data
 $u_0\in L^2_0(X,\nu)$ with $$\displaystyle \int_Xj^*_\gamma(u_0)d \nu<+\infty,$$  $T>0$, and any $f\in L^2(0,T;H^{-1}_m(X,\nu))$,  there exists a unique strong solution of the abstract Cauchy problem
\begin{equation}\label{ACP1}
 \left\{ \begin{array}{ll} u'(t) + \partial \Psi_\gamma(u(t))\ni f(t), \quad 0\le t\le T, \\[12pt] u(0) = u_0. \end{array} \right.
 \end{equation}
 \end{theorem}

 By the characterization of the operator $\partial_{H_m^{-1}(X,\nu)} \Psi_\gamma$ obtained in Theorem \ref{Charact1}, we have the unique strong solution  of Problem~\eqref{ACP1} is the only function $u \in C(0, T; H^{-1}(X,\nu))$ that satisfies
 \begin{equation}\label{PME12}
 \left\{ \begin{array}{lll}\frac{\partial u}{\partial t} - \Delta_m v(t) \ni f(t) \quad &\hbox{a.e} \ \ t \in (0, T), \\ \\ u(t) \in  \gamma(v(t)) &\hbox{a.e} \ \ t \in (0, T),\\ \\ u(0) = u_0.\end{array} \right.
 \end{equation}
 Therefore, the mild solutions given in Proposition~\ref{accretive} are in fact   strong solutions under the conditions of the above result. Under   this point of view we get moreover the following contraction principle: for $u_0,\widetilde{u_0}\in L^2_0(X,\nu)$ with $\displaystyle \int_Xj^*_\gamma(u_0)d \nu, \int_Xj^*_\gamma(\widetilde{u_0})d \nu<+\infty,$ and $f,\widetilde{f}\in L^1(0,T;H^{-1}_m(\Omega,\nu))$, and  $u,\widetilde{u}$ the respective strong solutions of  Problem~\eqref{PME12}, then
$$ ||u(t)-\widetilde u(t)||_{H^{-1}_m}\le ||u_0-\widetilde u_0||_{H^{-1}_m}
  +\int_0^t||f(s)-\widetilde f(s)||_{H^{-1}_m}ds,\quad\forall 0\le t\le T.
$$

\begin{remark}\label{1142}\rm Observe, that we can obtain the same result for data $u_0\in L^2(X,\nu)$ with $$\displaystyle \int_Xj^*_\gamma(u_0)d \nu<+\infty,$$    by using a translation argument. Set $$\displaystyle b=\frac{1}{\nu(X)}\int_Xu_0d\nu.$$
Since $j_\gamma^*$ is convex, we have
$$j_\gamma^*(b)=j_\gamma^*\left(\frac{1}{\nu(X)}\int_Xu_0d\nu\right)\le \frac{1}{\nu(X)}\int_Xj_\gamma^*(u_0)d\nu<+\infty.
$$
Set     $\widetilde\gamma$ defined via
   $$\widetilde{\gamma}^{-1}(s):= \gamma^{-1}(s+b) - (\gamma^{-1})^0(b).$$
   Then, for $\widetilde{u}_0:= u_0 - b$,
   $$\displaystyle \int_Xj^*_{\widetilde\gamma}(\widetilde u_0)d \nu= \int_Xj^*_\gamma(u_0)d \nu-\nu(X)j_\gamma^*(b)<+\infty.$$
Therefore, since  $$\widetilde{u}_0\in H^{-1}(X,\nu),$$ for $T>0$  and $f\in L^2(0,T;H^{-1}_m(X,\nu))$, there exists a unique strong solution $\widetilde{u}(t)$ of  problem
  \begin{equation}\label{ACP1New}
 \left\{ \begin{array}{ll} \widetilde{u}'(t) + \partial \Psi_{\widetilde{\gamma}}(\widetilde{u}(t))\ni f(t), \quad 0\le t\le T, \\[12pt] \widetilde{u}(0) = \widetilde{u}_0. \end{array} \right.
 \end{equation}
Then, there exists $\widetilde{v}(t)  \in \widetilde{\gamma}^{-1}(\widetilde{u}(t))$ such that
 $$\frac{\partial \tilde{u}}{\partial t} - \Delta_m \tilde{v}(t) \ni f(t) \quad \hbox{a.e} \ \ t \in (0, T).$$
 Now, if we define $u(t):= \widetilde{u}(t) + b$ and  $v(t):=\widetilde v(t)+ (\gamma^{-1})^0(b)$,  we have
 $$ v(t) =\widetilde v(t)+ (\gamma^{-1})^0(b) \in \widetilde{\gamma}^{-1}(\widetilde{u}(t)) + (\gamma^{-1})^0(b) = \gamma^{-1}(u(t)).$$
 Consequently, $u(t)$ is a strong solution of problem \eqref{PME12}.  $\blacksquare$
\end{remark}

\section{The Cahn-Hilliard Equations  on Random Walk Spaces}\label{Main}

\subsection{The Cauchy problem in $L^1$}\label{Cauchy}

  Let $[X,\mathcal{B},m^1,\nu_1]$  be a $m^1$-connected random walk space such that $\nu_1$ is reversible and  $0<\nu_1(X)<+\infty$.
 Let $[X,\mathcal{B},m^2,\nu_2]$  be a random walk space such that  $\nu_2$  is invariant and $0<\nu_2(X)<+\infty$. Assume moreover that
$$\nu_1\ll \nu_2$$
and $${\tt{R}}:=\frac{d\nu_1}{d\nu_2}\in L^\infty(X,\nu_2).$$
The above hypothesis implies that $L^1(X,\nu_2)\subset L^1(X,\nu_1)$ and
$${\tt m}|| f||_{L^1(X, \nu_1)}\le || f||_{L^1(X, \nu_2)},$$
with ${\tt m}=\frac{1}{||{\tt R}||_{L^\infty(X, \nu_2)}}.$ We moreover assume  that $L^1(X,\nu_1)$ is continuously embedded in  $L^1(X,{\nu_2})$ with
$$|| f||_{L^1(X, \nu_2)}\le {\tt M}|| f||_{L^1(X, \nu_1)}.$$

 We consider the doubly nonlocal system
\begin{equation}\label{DNLS1}
\left\{\begin{array}{lll} u_t(t,x) = \Delta_{m^1}  \mu(t,x) , \quad &(t,x) \in (0, \infty) \times X, \\[10pt] \mu(t,x)    \in  - \Delta_{m^2}u(t,x) + \partial F(u(t,x)), \quad &(t,x) \in (0, \infty) \times X, \\[10pt] u(0,x) = u_0(x), \quad &x \in  X,\end{array} \right.
\end{equation}
where
\begin{equation}\label{1137}
\partial F(r)  = \gamma^{-1}(r) - cr,
\end{equation}
being $\gamma$ a maximal monotone graph in $\R$ with $0\in\gamma(0)$ and $\gamma^-<\gamma^+$,   and $c>0$.

\begin{definition}{\rm   Given the random walks $m^1$, $m^2$, we define its {\it convolution} $m^1\ast m^2$ as the random walk defined by
$$\int_X \varphi (y) d(m^1\ast m^2)_x(y) := \int_X \left( \int_X \varphi (y) dm^2_z(y) \right) dm^1_x(z).$$
}
\end{definition}

\begin{lemma}\label{compp}       Let $ u \in  L^1(X,\nu_1) $.  We have
\begin{equation}\label{dobleLapl}
\Delta_{m^1}(\Delta_{m^2}u)(x) = \Delta_{m^1\ast m^2} u(x)- \Delta_{m^1}u(x) - \Delta_{m^2}u(x)  \quad \hbox{for all} \  x \in X.
\end{equation}
\end{lemma}
\begin{proof} We have
$$\Delta_{m^1}(\Delta_{m^2})u(x) = M_{m^1}(\Delta_{m^2}u)(x) - \Delta_{m^2}u(x) =  M_{m^1}(M_{m^2}u - u)(x) -  \Delta_{m^2}u(x).  $$
Now
$$M_{m^1}(M_{m^2}u)(x) = \int_X M_{m^2}u(z) dm^1_x(z)  = \int_X \left( \int_X u(y) dm^2_z(y) \right)dm^1_x(z) $$ $$= \int_X u(y) d(m^1\ast m^2)_x(y) = M_{m^1\ast m^2}u(x).$$
Hence
$$\Delta_{m^1}(\Delta_{m^2})u(x) =  M_{m^1\ast m^2}u(x) - M_{m^1}u(x) -  \Delta_{m^2}u(x),$$  that is
\[ \Delta_{m^1}(\Delta_{m^2})u(x) =  \Delta_{m^1\ast m^2}u(x) u(x) - \Delta_{m^1}u(x) - \Delta_{m^2}u(x).\hfill\qedhere
\]
\end{proof}
 As consequence of the above lemma, we can rewrite \eqref{DNLS1} as
 \begin{equation}\label{DNLS1rew}
\left\{\begin{array}{ll} u_t = \Delta_{m^1} v  -\Delta_{m^1\ast m^2} u + (1-c)\Delta_{m^1}u + \Delta_{m^2}u &\hbox{in } (0, \infty) \times X,
\\[10pt]
u \in \gamma(v) &\hbox{in } (0, \infty) \times X,
 \\[10pt] u(0)= u_0 &\hbox{in }  X.\end{array} \right.
\end{equation}

We define the operator $\mathbf{G}: L^1(X, \nu_1) \rightarrow L^1(X, \nu_1)$ as
 $$ \mathbf{G}(u)= \Delta_{m^1\ast m^2}u  + (c-1) \Delta_{m^1}u -\Delta_{m^2}u.$$
 As consequence of Remark \ref{remark.laplacian}, we have the following result.
 \begin{lemma}\label{Lipschitz} The operator $\mathbf{G}$ is Lipschitz continuous  in $L^1(X, \nu_1)$  and in $L^2(X, \nu_1)$.
 \end{lemma}

\begin{proof}
 By  Remark \ref{remark.laplacian} we have:
\\
1. $\Delta_{m^1}$ is 2-Lipschitz continuous   in $L^1(X, \nu_1)$  and in $L^2(X, \nu_1)$.
\\
2. $\Delta_{m^2}$ is 2-Lipschitz continuous   in $L^1(X, \nu_2)$  and in $L^2(X, \nu_2)$. Hence, for $f\in L^1(X,\nu_1)$,
$$
\displaystyle ||\Delta_{m^2}f||_{L^1(X, \nu_1)}\le \frac{1}{{\tt m}}||\Delta_{m^2}f||_{L^1(X, \nu_2)}
\le \frac{2}{{\tt m}}||f||_{L^1(X, \nu_2)}\le
2\frac{{\tt M}}{{\tt m}}||f||_{L^1(X, \nu_1)}.
$$
And, for $f\in L^2(X,\nu_1)$,  with a  similar argument,
$$\begin{array}{c}
\displaystyle ||\Delta_{m^2}f||_{L^2(X, \nu_1)}^2 \le 4\frac{{\tt M}}{{\tt m}}||f||_{L^2(X, \nu_1)}^2.
\end{array}
$$
3. For $f\in L^1(X,\nu_1)$,
$$
\displaystyle ||\Delta_{m^1}(\Delta_{m^2}f)||_{L^1(X, \nu_1)}\le 2  || \Delta_{m^2}f||_{L^1(X, \nu_1)}\le 4\frac{{\tt M}}{{\tt m}}||f||_{L^1(X, \nu_1)}.
$$
And, for $f\in L^2(X,\nu_1)$,
\[
\displaystyle ||\Delta_{m^1}(\Delta_{m^2}f)||_{L^2(X, \nu_1)}^2\le 4  || \Delta_{m^2}f||_{L^2(X, \nu_1)}^2\le 4^2\frac{{\tt M}}{{\tt m}}||f||_{L^2(X, \nu_1)}.\qedhere
\]
\end{proof}

Now, by means of the operator $${\mathbf B}_\gamma:={\mathbf B}^{m^1}_\gamma,$$ as   given in Definition \ref{operator1},  but for the random walk $[X,\mathcal{B},m^1,\nu_1]$,  we can rewrite \eqref{DNLS1rew} as the abstract Cauchy problem
\begin{equation}\label{ACProb}
\left\{ \begin{array}{ll} u'(t) + ({\mathbf B}_\gamma +  \mathbf{G})(u(t)) \ni 0,\quad t>0, \\[8pt] u(0) = u_0.
\end{array} \right.
\end{equation}
By Proposition \ref{accretive} and Lemma \ref{Lipschitz}, we have that $ {\mathbf B}_\gamma + \mathbf{G}+ L_{\mathbf G} I$ is an accretive operator   in $L^1(X, \nu_1)$, being $L_{\mathbf G}$ the Lipschitz constant of ${\mathbf G}$, with $\overline{D({\mathbf B}_\gamma + \mathbf{G}+ L_{\mathbf G} I)}^{ L^1(X, \nu_1)}=\left\{u\in L^1(X, \nu_1): \gamma^-\le u \le \gamma^+ \right\}$. We are going to see that Problem~\eqref{DNLS1rew}, via its abstract formulation given by~\eqref{ACProb}, has mild solutions for a large class of general initial data:

\begin{theorem}\label{1313}   Assume  that the random walk space  $[X,\mathcal{B},m^1,\nu_1]$  satisfies a {\it Poincar\'{e} inequality}.
  For  $u_0\in    L^1(X,\nu_1)$, $ \gamma^-\le u_0 \le \gamma^+ $,
such that
\begin{equation}\label{1257}\nu(X) \gamma^- < \int_X u_0 d\nu_1  < \nu(X) \gamma^+,
\end{equation}
Problem~\eqref{DNLS1rew} has a unique mild solution.
\end{theorem}

\begin{proof}
Let $u_0\in \left\{u\in L^1(X, \nu_1): \gamma^-\le u \le \gamma^+ \right\}$ satisfying condition~\eqref{1257}, and $0<\widetilde T<1/L_\mathbf{G}$, where $L_\mathbf{G}$ is the Lipstichz constant of $\mathbf{G}$ in $ L^1(X, \nu_1)$. And set, for a fixed $u_0$ as in the hypothesis, the functional
$$\mathfrak{F}:C([0,\widetilde T]: L^1(X, \nu_1)))\to C([0,\widetilde T]:L^1(X, \nu_1))$$ given by
$\mathfrak{F}(z)=u_z$ the mild solution of Problem~\eqref{ACProbPM} (or Problem~\eqref{PME2}) with initial datum $u_0$ and  $f=\mathbf{G}(z)$ (observe that $\mathbf{G}(z)\in L^1(0,\widetilde T; L^1(X, \nu_1))$. This functional is well defined thanks to Proposition~\ref{e1RCond1}, and, by item {\it 4.}  in such proposition and Lemma~\ref{Lipschitz}, it is contractive with constant $\widetilde TL_\mathbf{G}$, that we are assuming less than $1$. Then it has a unique fix point $u$ which  is the unique mild solution of Problem~\eqref{DNLS1rew}.
Indeed, for $\epsilon>0$, there exists $u_\epsilon$ a solution of an $\frac{\epsilon}{2(1+\widetilde TL_\mathbf{G})}-$discretization in $[0,\widetilde T]$ of
$ u' +  {\mathbf B}_\gamma u  \ni f$,
for $f=-\mathbf{G}(u)$, with $u_\epsilon(0)=u_0$. But since $\mathbf{G}$ is $L_{\mathbf{G}}$-Lipschitz continuous, this $u_\epsilon$ is a solution of an $\epsilon-$discretization in $[0,\widetilde T]$ of
$ u' +  ({\mathbf B}_\gamma+\mathbf{G}) u   \ni g$,
for $g=0$, with $u_\epsilon(0)=u_0$, and we have existence and uniqueness of mild solution   since  $ {\mathbf B}_\gamma + \mathbf{G}+ L_{\mathbf G} I$ is an accretive operator in $L^1(X, \nu_1)$ and the initial datum is in $\overline{D({\mathbf B}_\gamma + \mathbf{G}+ L_{\mathbf G} I)}^{L^1(X, \nu_1)}$.
Finally we can  extend the solutions up to any $T>0$ since continuity extension in time holds true for mild solutions.
\end{proof}

\begin{remark}\rm
We have that mild solutions are strong solutions when $\gamma^-$ and $\gamma^+$ are finite
 and
\begin{equation}\label{1258ini} \int_Xj_\gamma^*(u_0)d\nu_1<+\infty.
\end{equation}
Indeed,   since the mild solution $u$ is bounded between $\gamma^-$ and~$\gamma^+$, then we have that
  $\mathbf{G}(u)\in L^2(0,T; L^2(X, \nu_1))$. Therefore, we can apply {\it 5.} in Proposition~\ref{accretive}.

 In Subsection~\ref{lasecgf},  for $\nu_1=\nu_2=\nu$, we  see that   we have strong solutions  for any initial data in $ L^2(X, \nu)$ satisfying~\eqref{1258ini}.
$\blacksquare$
\end{remark}

\begin{remark}\rm
Observe that there is   mass preservation:
$$\int_\Omega u(t)d\nu_1=\int_\Omega u_0d\nu_1\quad\forall t\in[0,T].$$
 In fact, this is clear for strong solutions, but also for mild solutions, since this property is inherited from the stationary schemes.  $\blacksquare$
\end{remark}

\begin{example}\label{excon3}\rm

1. Consider the potential  $F_3(u)=  A(1 -u^2) + I_{[-1,1]}(u)$ given in the Introduction. Let $\gamma$ be  the inverse graph of $\partial I_{[-1,1]}$,
$$\gamma(r)=\left\{\begin{array}{ll} -1(=\gamma^-)&\hbox{if } r<0,\\[8pt]
[-1,1]&\hbox{if } r=0,\\[8pt]
1(=\gamma^+)&\hbox{if } r>0.
\end{array}\right.$$
Any initial datum $u_0 \in  L^1(X, \nu_1)$, $ -1\le u_0 \le 1$,
is in $L^2(X, \nu_1)$ and satisfies
\begin{equation}\label{1258}\int_Xj_\gamma^*(u_0)d\nu_1<+\infty,
\end{equation}
 since such condition, in this case, is equivalent to have $u_0(x)\in[-1,1]$, $x\in X$.   Therefore we will always have strong solutions for  any initial datum $u_0 \in  L^1(X, \nu_1)$, $ -1\le u_0 \le 1$, satisfying
$$-\nu_1(X)    < \int_X u_0 d\nu_1  < \nu_1(X).$$
But, if $u_0$
satisfies  $$ -\nu_1(X) = \int_X u_0 d\nu_1,$$
which is equivalent to say that the initial datum is the pure phase $$u_0=-1,$$ then the pure phase $u(t)=-1$ is the solution of~Problem~\eqref{DNLS1rew}. Similarly, for the another pure phase $u_0=1$, $u(t)=1$ is the solution of the problem.

2.   For the case with the potential $F_1(r) = ((1+r) \log (1+r)+ (1-r)  \log (1-r)) - \frac{c}{2}r^2$,  and datum  $u_0 \in  L^1(X, \nu_1)$, $ -1\le u_0 \le 1$,
$$ -\nu_1(X)    < \int_X u_0 d\nu_1  < \nu_1(X),$$
and
$$\int_Xj_\gamma^*(u_0)d\nu_1<+\infty,
$$
where here $\gamma^{-1}(r)=\log(1+r) - \log(1-r)$, we will  have strong solutions.

3. Consider now the regular polynomial potential $F_2(u) = \frac14(u^2 -1)^2$. In this case,   $\gamma^{-1}(r)=r^3$, and  $\gamma^-=-\infty$ and $\gamma^+=+\infty$. Therefore we have existence of mild solutions for any initial datum   in $L^1(X, \nu_1)$.
For $\nu_1=\nu_2=\nu$ and for data in $ L^2(X, \nu)$ satisfying $ \int_Xj_\gamma^*(u_0)d\nu<+\infty$, which in this case is equivalent to ask for
 data   $u_0\in L^4(X, \nu)$, we also have strong solutions on account of the results given in Subsection~\ref{lasecgf}.  $\blacksquare$
\end{example}

\subsection{The regional Neumann problem}
 Let $[X,\mathcal{B},m^1,\nu_1]$  be a   random walk space such that $\nu_1$ is reversible and $\nu_1(X)<+\infty$.
 Let $[X,\mathcal{B},m^2,\nu_2]$  be a random walk space such that $\nu_1$ is invariant and $\nu_2(X)<+\infty$.

 Let  $\Omega \in \mathcal{B}$ with $0<\nu_i(\Omega)<+\infty$, $i=1,2$. Consider  the random walk spaces
$$  \left[\Omega,\mathcal{B}_\Omega,(m^i)^{\Omega},  \nu_i \res \Omega \right], \quad i=1,2,$$
given in Example \ref{example.restriction.to.Omega}; assume $\left[\Omega,\mathcal{B}_\Omega,(m^1)^{\Omega},  \nu_1\res \Omega \right] $ is $(m^1)^\Omega$-connected.
  Assume moreover that
$$\nu_1\res \Omega\ll \nu_2\res \Omega,$$
and $$\frac{d\nu_1\res \Omega}{d\nu_2\res \Omega}\in L^\infty(\Omega,\nu_2).$$
Assume also that  that $L^1(\Omega,\nu_1)$ is continuously embedded in  $L^1(\Omega,\nu_2)$.

Then, if we apply  the results of Subsection \ref{Cauchy} to these random walk spaces, we have that Problem \eqref{DNLS1} corresponds to the following Cahn-Hilliard problem
\begin{equation}\label{DNLS1N}
\left\{\begin{array}{lll} u_t(t,x) = \displaystyle\int_\Omega ( \mu(t,y) - \mu(t,x)) d(m^1)_x(y), \  &(t,x) \in (0, \infty) \times \Omega, \\[10pt] \mu(t,x)  \in - \displaystyle\int_\Omega (u(t,y)-u(t,x)  d(m^2)_x(y) + \partial F(u(t,x)), \  &(t,x) \in (0, \infty) \times \Omega, \\[10pt] u(0,x) = u_0(x), \quad &x \in  \Omega.\end{array} \right.
\end{equation}
As consequence of Theorem \ref{1313}    and Theorem~\ref{ExistUniqStrong}  below, we have the following result about existence and uniqueness of solutions to Problem \eqref{DNLS1N}.

\begin{theorem}\label{Wellpost}    Assume, joint to the previous assumptions,  that the random walk space $ \left[\Omega,\mathcal{B}_\Omega,(m^1)^{\Omega},  \nu_1 \res \Omega \right]$  satisfies a {\it Poincar\'{e} inequality}.
  For  $u_0\in    L^1(\Omega,\nu_1)$, $ \gamma^-\le u_0 \le \gamma^+$,
  such that
\begin{equation}\label{1257N} \nu(\Omega) \gamma^- < \int_\Omega u_0 d\nu_1  < \nu(\Omega) \gamma^+,
\end{equation}
Problem~\eqref{DNLS1N} has a unique mild solution.

 Moreover,  if  $\nu_1\res \Omega=\nu_2\res \Omega$ and it is also reversible with respect to $m^2$,
for $u_0\in     L^2(\Omega, \nu_1)$
such that
\begin{equation}\label{1258B}  \int_\Omega j_\gamma^*(u_0)d\nu_1<+\infty,
\end{equation}
Problem~\eqref{DNLS1N} has a unique strong solution.
\end{theorem}


\subsection{The gradient flow}\label{lasecgf}   Fife showed  in \cite{Fife}  that the local system~\eqref{Ch1} is the gradient flow of the Ginzburg-Landau free energy
$$\mathcal{E}_l (u) = \frac{1}{2} \int_\Omega \vert \nabla u \vert^2 dx + \int_\Omega F(u) dx,$$
respect the the scalar product
$$\langle v_1, v_2 \rangle_{H^{-1}} := \langle \nabla \phi_1, \nabla \phi_2 \rangle_{L^2},$$
where  $\displaystyle v_1, v_2\in H^{-1}:= \left\{ v \in L^2(\Omega)  :  \int_\Omega v =0 \right\}$ and $\phi_1, \phi_1 \in H^{-1}$ are the unique solutions of the Neumann problem
$$\left\{\begin{array}{ll} \Delta \phi_i = v_i \quad  &\hbox{in} \ \ \Omega, \\[10pt] \displaystyle\frac{ \partial  \phi_i}{\partial \eta} =0 \quad  &\hbox{on} \ \ \partial \Omega. \end{array} \right.$$

 This approach is frequently used in the literature to deal with the different versions (local-nonlocal, nonlocal-nonlocal) of the Cahn-Hilliard system. We will obtain a similar development in the ambient space of random walk spaces.

  Let $[X,\mathcal{B},m^1,\nu]$  be a $m^1$-connected random walk space such that $\nu$ is reversible with respect to $m^1$ and $\nu(X)<+\infty$.
 Let $[X,\mathcal{B},m^2,\nu]$  be a random walk space such that $\nu$ is  reversible with respect to $m^2$. Observe we are imposing the  restrictive assumption that both spaces have the  same invariant and reversible measure. This assumptions are always satisfied if we are dealing with nonlocal kernels $m^{J_1}$ and $m^{J_2}$ in $\mathbb{R}^N$ like in Example~\ref{example.nonlocalJ} since the Lebesgue measures is an invariant and reversible measure for both, they are also satisfied for certain   graphs with different, but related, weights  defining $m^1$ and $m^2$.

Since we are dealing with potentials $F$ such that $\partial F(r)=\gamma^{-1}(r)-cr$, with $\gamma^{-1}$ possibly multivalued, the Ginzburg-Landau free energy functional is defined as $$\mathcal{E} :   H_{m^1}^{-1}(X, \nu) \rightarrow ]-\infty, +\infty],$$ with
\begin{equation}\label{GLEnerg}
\mathcal{E}(u):= \mathcal{H}(u)+\Psi(u) - \frac{c}{2}\int_X u^2  d \nu,
\end{equation}
 where $\mathcal{H}:  H^{-1}_{m^1}(X, \nu)\rightarrow \mathbb{R}$ is
$$ \mathcal{H}(u):=  \frac{1}{4} \int_{X \times X} (u(x) - u(y))^2 dm^2_x(y) d\nu(x),$$
  $ \Psi :  H^{-1}_{m^1}(X, \nu)\rightarrow ]\infty, +\infty]$ is
\begin{equation}\label{1744adapted}\Psi(u):= \left\{ \begin{array}{ll} \displaystyle\int_X j^*_\gamma(u(x)) d\nu(x) \quad &\hbox{if} \ \   j^*_\gamma(u) \in L^1(X, \nu), \\ \\ + \infty \quad &\hbox{otherwise}.\end{array} \right.
\end{equation}
Set \begin{equation}\label{GLEnerg1}
\mathcal{K}(u):=  \frac{c}{2}\int_X u^2 d\nu.
\end{equation}

Applying Theorem \ref{Charact1}, we have the following result.

\begin{lemma}\label{lsw1}
The functional $\Psi$ is convex and lower semi-continuous in $ H_{m^1}^{-1}(X, \nu)$ and
\begin{equation}\label{e2ChactN}\begin{array}{ll}
\partial_{H_{m^1}^{-1}(X, \nu)} \Psi=  &\Big\{ (u, w) \in  H_{m^1}^{-1}(X, \nu) \times H_{m^1}^{-1}(X, \nu) :  w = - \Delta_{m^1} v, \\[8pt] &\qquad v \in   L^2(X, \nu), \ u(x) \in \gamma(v(x)) \   \nu\hbox{-a.e.} \ x \in X \Big\}.\end{array}
\end{equation}
\end{lemma}

\begin{lemma}\label{lsw2} The operator $\mathcal{H}$ is  proper, convex and continuous in $H_{m^1}^{-1}(X, \nu)$.
\end{lemma}
\begin{proof}   Obviously $\mathcal{H}$ is  proper and convex, and the continuity follows directly from  Proposition~\ref{isomorph}. In fact, let $\{u_n\}_n\subset  H_{m^1}^{-1}(X, \nu)$ and $u\in H_{m^1}^{-1}(X, \nu)$ be such that $$||.||_{H_{m^1}^{-1}(X, \nu)}-\lim_n u_n =u.$$ Then, from Proposition~\ref{isomorph},
$$u_n\to u\quad \hbox{in }   L^2(X,\nu).$$
Therefore, on account of this convergence, we get
\[\mathcal{H}(u) =\lim_n \mathcal{H}(u_n).\qedhere\]
\end{proof}

\begin{theorem}\label{subdiff1}
We have that $\mathcal{H}+\Psi$ is convex and lower semi-continuous in $H^{-1}_{m^1}(X, \nu)$. Moreover, $w \in \partial_{H_{m^1}^{-1}(X, \nu)} (\mathcal{H}+\Psi)(u)$ if and only if
\begin{equation}\label{e1subdiff1}
  \exists v \in  L^2(X, \nu), \hbox{with} \ u(x) \in \gamma(v(x))\ \hbox{$\nu$-a.e., such that} \ w = \Delta_{m^1}(\Delta_{m^2} u) - \Delta_{m^1} v.
\end{equation}
\end{theorem}
\begin{proof}
By Lemma \ref{lsw1} and Lemma \ref{lsw2}, we have   $\mathcal{H}+\Psi$ is convex and lower semi-continuous in $ H^{-1}_{m^1}(X, \nu)$. Now, since
$${\rm int}(D(\mathcal{H})) \cap D(\Psi) =    L^2_0(X,\nu)   \cap L^1(X, \nu) \not= \emptyset,$$
applying \cite[Corollary 2.11]{Brezis}, we have
$$ \partial_{H_{m^1}^{-1}(X, \nu)} (\mathcal{H}+\Psi) = \partial_{ H_{m^1}^{-1}(X, \nu)} \mathcal{H} + \partial_{ H_{m^1}^{-1}(X, \nu)} \Psi.$$
Hence, $w \in  \partial_{H_{m^1}^{-1}(X, \nu)} (\mathcal{H}+\Psi)(u)$ if and only if $$\exists w_1 \in \partial_{H_{m^1}^{-1}(X, \nu)} \mathcal{H}(u) \ \hbox{and} \ \exists   w_2 \in \partial_{H_{m^1}^{-1}(X, \nu)} \Psi(u) \ \hbox{such that} \ w = w_1 + w_2.$$

Now, by Lemma \ref{lsw1}, there exists $v \in L^2(X, \nu), \ u(x) \in \gamma(v(x)) \   \nu\hbox{-a.e.}  \ x \in X$, such that $w_2 = - \Delta_{m^1} v$.

On the other hand, since $\mathcal{H}$ is convex, we have
$
w_1 \in \partial_{ H_{m^1}^{-1}(X, \nu)} \mathcal{H}(u)$ in and only if $$\liminf_{t \to 0^+} \frac{1}{t} \left(\mathcal{H}(u + t h) - \mathcal{H}(u) \right)  \geq \langle w_1, h \rangle_{H_{m^1}^{-1}(X, \nu)} \quad \forall h \in H^{-1}_{m^1}(X, \nu). $$
 Now, for $h \in   H^{-1}_{m^1}(X, \nu)$, we have
$$ \liminf_{t \to 0^+} \frac{1}{t} \left(\mathcal{H}(u + t h) - \mathcal{H}(u) \right) =  \displaystyle\frac12 \int_{X \times X}  \nabla u(x,y)   \nabla h(x,y)  d(\nu\otimes (m^2)_x)(x,y) $$ $$ = -\int_X \Delta_{m^2} u h  d\nu = -\int_X \Delta_{m^1}^{-1}(\Delta_{m^1}(\Delta_{m^2} u)) h d\nu = \langle \Delta_{m^1}(\Delta_{m^2} u), h \rangle_{ H^{-1}_{m^1}(X, \nu)}.$$
Therefore, $w_1 = \Delta_{m^1}(\Delta_{m^2} u)$.
\end{proof}

As consequence of Theorem \ref{subdiff1}, we can  rewrite Problem~\eqref{DNLS1} as
\begin{equation}\label{GRAD1}
\left\{ \begin{array}{ll} u'(t) + \partial_{H_{m^1}^{-1}(X, \nu)} (\mathcal{H}+\Psi))(u(t)) - \mathcal{B}(u(t))\ni 0,\quad t>0, \\[8pt] u(0) = u_0.
\end{array} \right.
\end{equation}
being
$$ \mathcal{B}(u):=  - c \Delta_{m^1} u.$$

\begin{proposition}\label{operB} The operator $\mathcal{B}$ is $H^{-1}_{m^1}$-Lipschitz continuous.
\end{proposition}
\begin{proof} Since $\mathcal{B}$ is $L^2(X, \nu)$-Lipschitz continuous, the result is consequence of  Proposition~\ref{isomorph}.
\end{proof}

\begin{theorem}\label{ExistUniqStrong}  Assume  that the random walk space  $[X,\mathcal{B},m^1,\nu]$  satisfies a Poincar\'{e} inequality.
Then, for any  $u_0\in L^2(X, \nu)$
with
\begin{equation}\label{1258New}\int_Xj_\gamma^*(u_0)d\nu<+\infty,
\end{equation}
there exists a unique   strong solution of Problem \eqref{DNLS1}.
Moreover,  the mild solution  given in Theorem~\ref{1313} is such strong solution under the above conditions.
\end{theorem}

\begin{proof}
Having in mind Theorem \ref{subdiff1} and Proposition \ref{operB}, applying \cite[Proposition 3.12]{Brezis}, we have the   result for initial data
 in $L^2_0(X,\nu)$.
Now, using the translation argument  given in Remark~\ref{1142}, we also have  strong solutions for data
$u_0\in L^2(X, \nu)$ since we have  that $-c\Delta_{m^1}u\in L^2(0,T;H_{m^1}^{-1}(X, \nu))$.
\end{proof}

\begin{remark}\label{gdntflw}\rm   Problem~\eqref{DNLS1}, or its equivalent expression given in Problem~\eqref{GRAD1}, is the gradient flow in $H_{m^1}^{-1}(X, \nu)$ of the Ginzburg-Landau free energy functional $\mathcal{E}$ given in~\eqref{GLEnerg}:
\begin{equation}\label{GRAD2}
\left\{ \begin{array}{ll} u'(t) + \partial_{H_{m^1}^{-1}(X, \nu)}  \mathcal{E} (u(t))\ni 0,\quad t>0, \\[12pt] u(0) = u_0.
\end{array} \right.
\end{equation}
 Indeed, by definition (see~\eqref{subdif1}),
$h\in \partial_{H_{m^1}^{-1}(X, \nu)}  \mathcal{E}(u)$ if and only if
$$\liminf_{t \to 0^+} \frac{\mathcal{E}  (u + tv) - \mathcal{E}  (u)}{t} \geq \langle h, v \rangle_{H_{m^1}^{-1}(X, \nu)} \quad   \forall \, v \in H_{m^1}^{-1}(X, \nu).
$$
Now, by convexity,    $\widetilde h\in \partial_{H_{m^1}^{-1}(X, \nu)}(\mathcal{H}+\Psi) (u)$ if and only if
\begin{equation}\label{add01}\liminf_{t \to 0^+} \frac{(\mathcal{H}+\Psi) (u + tv) - (\mathcal{H}+\Psi) (u)}{t} \geq \langle \widetilde h, v \rangle_{H_{m^1}^{-1}(X, \nu)}  \quad \forall \, v \in H_{m^1}^{-1}(X, \nu).
\end{equation}
And, on the other hand, for $v \in H_{m^1}^{-1}(X, \nu)$,
\begin{equation}\label{add02}
\begin{array}{c}\displaystyle\exists\lim_{t \to 0^+} \frac{\displaystyle - \frac{c}{2}\int_X  (u + tv)^2d\nu  + \frac{c}{2}\int_X  u^2 d\nu}{t}=
-\langle  cu, v \rangle_{L^2(X, \nu)}\\ \\=-\langle  c\Delta_{m^1}^{-1}\Delta_{m^1}u, v \rangle_{L^2(X, \nu)}=\langle  c\Delta_{m^1}u, v\rangle_{H_{m^1}^{-1}(X, \nu)}.
\end{array}
\end{equation}
Therefore adding up~\eqref{add01} and~\eqref{add02}, \eqref{GRAD1} can be written as~\eqref{GRAD2}.  $\blacksquare$
\end{remark}

  \section{Some properties of the solution}\label{Properties}

In this section we obtain some properties of the solutions and of their asymptotic behaviour.

In the next result we see that, for $L^\infty(X,\nu)$-bounded initial data, strong solutions stay bounded in $L^\infty(0,T;L^\infty(X,\nu))$, therefore  the   solutions obtained in Theorem~\ref{Wellpost} coincide with  the solutions obtained, with a different method, by Gal and Shomberg in~\cite{GS} for the case considered in such paper. Concretely, they impose that the potential $F \in C^2(\R)$ satisfies $F'(0) =0$ and
$$F^{''} (r) + a_K(x) \geq  c_0  \quad \hbox{for all $r \in \R$ and a.e.   $x \in \Omega$,}$$
for some constant $c_0 >0$, where $\displaystyle a_K(x):= \int_\Omega K(x-y) dy$ (being $K$   the interaction kernel playing the role of $m^2$).
  This case includes, for example, the double-well potential $F_2$   given in the Introduction   but not $F_1$ or $F_3$.
Therefore the above existence and uniqueness result generalize the results by Gal and Shomberg. At our knowledge, Theorem~\ref{1313} and  Theorem~\ref{1258New} give new existence and uniqueness results for  the Cahn-Hilliard problem by the generality of the potentials and the transitions that can be considered on the system, and the large class of initial data.

\begin{proposition}\label{tento3}
  Under the conditions in Theorem~\ref{ExistUniqStrong}, assume moreover that  $u_0\in   L^p(X,\nu)$, $2\le p\le +\infty$. Then, for $u$ being the unique strong solution of Problem \eqref{DNLS1}, there exists a constant $C>0$ (which depends only on $c$ given in~\eqref{1137}) such that
 $$||u(t)||_{L^p(X,\nu)}\le ||u_0||_{L^p(X,\nu)}e^{CT}\quad\hbox{for all} \ \  0<t<T.$$
\end{proposition}

\begin{proof}
We know that $u\in W^{1,1}_{loc}((0,T);L^2(X,\nu)\cap C([0,T]:L^2(X,\nu)$, $u(0)=u_0$, and  there exists $v(t)\in L^2(X,\nu)$, $v(t)\in\gamma(u(t))$ $\nu_1$-a.e, such that, for almost every $t\in(0,T)$,
$$
  u_t = \Delta_{m^1} v  -\Delta_{m^1\ast m^2} u + (1-c)\Delta_{m^1}u + \Delta_{m^2}u \quad\hbox{in } (0, T) \times X.
$$
Let $p\geq 2$  and $j_{p,k}(r)$ be the primitive of the nondecreasing function $r|T_k(r)|^{p-2}$. Multiplying the above equation by
 $u|T_k(u)|^{p-2}$, integrating over $X$, and misleading nonnegative terms, we get
 \begin{equation}\label{blue01}\begin{array}{c}\displaystyle\frac{d}{dt}\int_X j_{k,p}(u(t))d\nu\le   \int_{X }    -\Delta_{m^1\ast m^2} u(t,x)\left(u|T_k(u)|^{p-2}\right)(t,x) d\nu(x)  \\
\\
\displaystyle + (c-1)^+  \int_{X}   -\Delta_{m^1} u(t)\left(u(t)|T_k(u(t))|^{p-2}\right)  d\nu
\\ \\
\displaystyle = \int_{X}\int_{X}    -(u(t,y)-u(t,x))\left(u(t,x)|T_k(u(t,x))|^{p-2}\right) d(m^1\ast m^2)_x(y)d\nu(x)  \\
\\
\displaystyle + (c-1)^+  \int_{X}\int_{X}   -(u(t,y)-u(t,x))u(t,x)\left(u(t,x)|T_k(u(t,x))|^{p-2}\right) dm^1_x d\nu(x)
\end{array}
\end{equation}
Now,
\begin{equation}\label{cabreo2}  -(a-b)b \vert T_k(b) \vert^{p-2}  \le (|a|+|b|)|b|\vert T_k(b) \vert^{p-2} \leq 2 a^2 \vert T_k(a) \vert^{p-2} + 2 b^2 \vert T_k(b) \vert^{p-2},
\end{equation}
and, a direct calculation shows that
\begin{equation}\label{cabreo1}r^2|T_k(r)|^{p-2}\le pj_{k,p}(r).
\end{equation}
Then, by \eqref{cabreo2} and \eqref{cabreo1},  \eqref{blue01}  yields
\begin{equation}\label{blue01dos}\begin{array}{c}\displaystyle\frac{d}{dt}\int_X j_{k,p}(u(t))d\nu\le 2p  \displaystyle  \int_{X}\int_{X}    \left(j_{k,p}(u (t,x))+j_{k,p}(u(t,y))\right) d(m^1\ast m^2)_x(y)d\nu(x)  \\
\\
\displaystyle + 2p (c-1)^+  \int_{X}\int_{X}   \left(j_{k,p}(u(t,x))+j_{k,p}(u(t,y))\right) dm^1_x(y) d\nu(x)
\end{array}
\end{equation}
 Now, it is easy to see that $\nu$ is invariant with respect to $m^1\ast m^2$.
Therefore, using the invariance of $\nu$ with respect to $m^1$ and $m^1\ast m^2$, \eqref{blue01dos} gives
$$\frac{d}{dt}\int_X j_{k,p}(u(t))d\nu\le Cp\int_X j_{k,p}(u(t))d\nu,$$
with $C=4\max\{c,1\}$. Then, from Gr\"onwall's lemma,
$$\int_X j_{k,p}(u(t))d\nu\le e^{CpT}\int_X j_{k,p}(u_0)d\nu \quad\forall 0\le t\le T.$$
Hence
$$\left(\int_X pj_{k,p}(u(t))d\nu\right)^{1/p}\le e^{CT}\left(\int_X pj_{k,p}(u_0)d\nu\right)^{1/p} \quad\forall 0\le t\le T,$$
and,  taking  limits as $k$ goes to $+\infty$ in the above expression we get
\begin{equation}\label{lpest}\left(\int_X |u(t)|^pd\nu\right)^{1/p}\le e^{CT}\left( \int_X |u_0|^pd\nu\right)^{1/p} \quad\forall 0\le t\le T.
\end{equation}
This gives the proof if  $p$ is finite. For $p=+\infty$,
 we have just to  take limits as $p$ goes to $+\infty$ in~\eqref{lpest} to get
\[||u(t)||_{L^\infty(X,\nu)}\le ||u_0||_{L^\infty(X,\nu)}e^{CT}\quad\forall 0<t<T.\qedhere\]
\end{proof}

\begin{remark}{\rm
  A similar result can be obtained for strong solutions of Problem \eqref{DNLS1} under the general conditions of Subsection~\ref{Cauchy}.}  $\blacksquare$
\end{remark}

  Let us see   that the  solution to Problem~\eqref{DNLS1} satisfies an energy identity. This will allow to get an uniform estimate in time for the $L^2$-norm of the solution    under natural extra conditions (see Corollary~\ref{ACOT11}).

\begin{proposition}\label{sufn} Under the conditions in Theorem~\ref{ExistUniqStrong}, for $u$ the strong solution to Problem~\eqref{DNLS1}, set
$\widetilde{\mathcal{E}} (t):= \mathcal{E}(u(t))$. Then
 $$\frac{d}{dt}\widetilde{\mathcal{E}} (t)=-\frac {1}{2} \int_{X \times X} |\nabla \mu(t)|^2 d(\nu\otimes (m^1)_x) \quad\hbox{for a.e. } t>0,$$
 where
 $\mu    =  - \Delta_{m^2}u + v -cu $, $v\in  \gamma^{-1}(u)$,  as corresponding to the definition of strong solution.
\end{proposition}

\begin{proof}
 We have that
\begin{equation}\label{GRAD1wh}
\left\{ \begin{array}{ll} u'(t) + \partial_{H_{m^1}^{-1}(X, \nu)} (\mathcal{H}+\Psi))(u(t)) - c \Delta_{m^1} u \ni 0,\quad t>0, \\[8pt] u(0) = u_0.
\end{array} \right.
\end{equation}
Now, by \cite[Lemme 3.3]{Brezis}, for almost every $t>0$,
$$\frac{d}{dt}\Big((\mathcal{H}+\Psi))(u(t))\Big)=\langle h,u'(t)\rangle_{H_{m^1}^{-1}(X, \nu)}\quad\forall h\in \partial_{H_{m^1}^{-1}(X, \nu)} (\mathcal{H}+\Psi))(u(t)).$$
Then, since $$ c \Delta_{m^1} u-u'(t)\in \partial_{H_{m^1}^{-1}(X, \nu)} (\mathcal{H}+\Psi))(u(t))$$ and
$$  u'(t)=\Delta_{m^1}\mu(t),$$
we have
$$\begin{array}{c}\displaystyle
\frac{d}{dt}\Big((\mathcal{H}+\Psi)(u(t))\Big)=\langle c \Delta_{m^1} u-u'(t),u'(t)\rangle_{H_{m^1}^{-1}(X, \nu)}
\\ \\\displaystyle
=\langle c \Delta_{m^1} u,u'(t)\rangle_{H_{m^1}^{-1}(X, \nu)}  -\langle \Delta_{m^1}\mu(t),\Delta_{m^1}\mu(t)\rangle_{H_{m^1}^{-1}(X, \nu)}
\\ \\\displaystyle
=\frac{d}{dt}\Big(\mathcal{K}(u(t)) \Big) -\langle \Delta_{m^1}\mu(t),\Delta_{m^1}\mu(t)\rangle_{H_{m^1}^{-1}(X, \nu)},
\end{array}
$$
and consequently,
$$\begin{array}{c}\displaystyle\frac{d}{dt}\Big((\mathcal{H}+\Psi -\mathcal{K})(u(t))\Big) = -\langle \Delta_{m^1}\mu(t),\Delta_{m^1}\mu(t)\rangle_{H_{m^1}^{-1}(X, \nu)},
\\ \\
\displaystyle
=-\frac {1}{2} \int_{X \times X} |\nabla \mu(t)|^2 d(\nu\otimes (m^1)_x).
\end{array}
$$
\end{proof}

\begin{corollary}\label{ACOT11}    Assume the conditions in Theorem~\ref{ExistUniqStrong} and suppose  $[X,\mathcal{B},m^2,\nu]$   satisfies a Poincar\'{e} inequality and the potential $\displaystyle j_\gamma^*(r)-\frac{c}{2}r^2 $ is bounded from bellow. For $u$ the strong solution to Problem~\eqref{DNLS1},
we have that
$$\{u(t):t>0\} \hbox{ is bounded in $L^2(X,\nu)$.}$$
\end{corollary}
\begin{proof}   As a consequence of Proposition~\ref{sufn} we have that   $\widetilde{\mathcal{E}} (t)$  is absolutely continuous and nonincreasing, therefore
$$\mathcal{E}(u(t))\le \mathcal{E} (u_0)\quad\forall t>0.$$
Hence, since the potential $F$ is bounded from below, we have
 \begin{equation}\label{1814}\mathcal{H}(u(t)) \leq C, \quad \hbox{for all $t \geq0$}.\end{equation}
 Finally, using the Poincar\'{e} inequality for $[X,\mathcal{B},m^2,\nu]$, from~\eqref{1814}, we get the thesis.
\end{proof}

\subsection{Asymptotic behaviour}

 In~\cite{GS} a very interesting analysis is done on the asymptotic behaviour of the solutions  for  the case of smooth convolution kernels and smooth potential. Their assumptions guaranties a regularization  and  a \L{}ojasiewicz-Simon inequality. With the generality of the nonlocal interactions and potentials studied here it is not clear if one can  get such tools.   Nevertheless we can prove some facts.

From now on we will assume that we are under the conditions in Theorem~\ref{ExistUniqStrong}.
 Then, we can define a semigroup $T(t): L^2(X, \nu) \rightarrow L^2(X, \nu)$, such that, for every $u_0 \in L^2(X, \nu)$, $T(t)u_0:= u(t)$ is the unique strong solution of problem~\eqref{DNLS1}.  For $u_0 \in L^2(X, \nu)$, we define its {\it omega limit set}
$$\omega(u_0):=   \{ w \in L^2(X, \nu)   : \exists t_n \to +\infty \ \hbox{such that} \ T(t_n) u_0 \to w  \hbox{ in $L^2(X,\nu)$} \}.$$

We also consider  the set of {\it equilibria solutions} of the problem,
$$\mathbb{E}=\mathbb{E}(T(t)):= \{ u \in L^2(X, \nu) \  : \ T(t)u = u \quad \forall \, t \geq 0 \} $$
 that can be characterized as $$  \mathbb{E}=
\{ u \in L^2(X, \nu)   : \exists  \mu  \hbox{ constant   with } \mu+\Delta_{m^2}u +cu  \in \gamma^{-1}(u) \}.$$

\begin{theorem}\label{ASINTB}
Assume  that the random walk space  $[X,\mathcal{B},m^2,\nu]$  also satisfies a Poincar\'{e} inequality.
Let $u_0 \in L^2(X, \nu)$ be and assume that
\begin{equation}\label{CCOnd2}  \gamma^-     < \frac{1}{\nu(X)} \int_X u_0 d\nu  <  \gamma^+,
\end{equation}
and also that the  omega limit set $\omega(u_0)$ is not empty. Then,
\begin{equation}\label{gradlike}
\omega(u_0) \subset \mathbb{E}.
\end{equation}
\end{theorem}
\begin{proof} Let $u(t):= T(t)u_0$ be,  and $\mu    =  - \Delta_{m^2}u + v -cu $, $v\in  \gamma^{-1}(u)$,  as corresponding to the definition of strong solution.

Given $u_\infty\in \omega(u_0)$, there exists  a sequence  $t_n \to +\infty$  such that
$$u(t_n)\to u_\infty \hbox{ in $L^2(X,\nu)$}.$$
From   Proposition~\ref{sufn} and the Poincar\'{e} inequality  we have that
\begin{equation}\label{forco1}\int_{0}^\infty\int_X|\mu(s,x)-\overline{\mu}(s)|^2 d\nu(x)ds\le C,
\end{equation}
where $\displaystyle \overline{\mu}(s)=\frac{1}{\nu(X)}\int_X\mu(s,x)d\nu(x)$,
and
\begin{equation}\label{forco2}\alpha_n:=\int_{t_n}^\infty\int_X|\mu(s,x)-\overline{\mu}(s)|^2 d\nu(x)ds\le  2\lambda_{m^2}(\widetilde{\mathcal{E}}(t_n)-\widetilde{\mathcal{E}}_\infty)\to 0,
\end{equation}
where
$\widetilde{\mathcal{E}}_\infty:= \lim_{t\to +\infty}\widetilde{\mathcal{E}} (t).$

 Let us see that there exists $\widetilde t_n\in[t_n,t_n+C/\sqrt{\alpha_n}]$ such that
\begin{equation}\label{pere01}\int_X|\mu(\widetilde t_n,x)-\overline{\mu}(\widetilde t_n)|^2 d\nu(x)\le \sqrt{\alpha_n}.
\end{equation} Indeed, arguing by contradiction, if for almost all $s\in [t_n,t_n+C/\sqrt{\alpha_n}]$ we have
$$\int_X|\mu(s,x)-\overline{\mu}(s)|^2 d\nu(x)> \sqrt{\alpha_n},
$$
then,
$$\int_{t_n}^{t_n+C/\sqrt{\alpha_n}}\int_X|\mu(s,x)-\overline{\mu}(s)|^2 d\nu(x)> C,
$$
which gives a contradiction with~\eqref{forco1}.

  Let us now see that
\begin{equation}\label{parex01}u(\widetilde t_n)\to u_\infty \hbox{ in $L^2(X, \nu)$}.
\end{equation}
In fact,
$$||u(\widetilde t_n)-u(t_n)||_{L^2(X,\nu)}=\left\|\int_{t_n}^{\widetilde t_n}{ \frac{d}{ds}}u(s)ds\right\|_{L^2(X,\nu)}$$
$$=\left\|\int_{t_n}^{\widetilde t_n}\Delta_{m_1}\mu(s)  ds\right\|_{L^2(X,\nu)}=
\left\|\int_{t_n}^{\widetilde t_n} \Delta_{m_1}(\mu(s)-\overline{\mu}(s))  ds\right\|_{L^2(X,\nu)}$$
$$\le  (t_n-\widetilde t_n)^{1/2}\left(\int_{t_n}^{\widetilde t_n}  \left\| \Delta_{m_1}(\mu(s)-\overline{\mu}(s))\right\|_{L^2(X,\nu)}^2   ds\right)^{1/2}$$
$$\le 2 (t_n-\widetilde t_n)^{1/2}\left(\int_{t_n}^{+\infty} \left\| \mu(s)-\overline{\mu}(s)\right\|_{L^2(X,\nu)}^2   ds\right)^{1/2}$$
$$=2\ (t_n-\widetilde t_n)^{1/2}\left(\alpha_n\right)^{1/2}\le 2 C^{1/2}\left(\alpha_n\right)^{1/4} \to 0.$$

Now we see that $\{\overline{\mu}(\widetilde t_n)\}_n$ is bounded. In fact, if
$$\lim \overline{\mu}(\widetilde t_n)=+\infty,$$ then, from~\eqref{pere01},  we have that
$$  \mu(\widetilde t_n,.)  \to +\infty\quad\hbox{$\nu$-a.e.}$$ But, since
$$\mu(\widetilde t_n)    =  - \Delta_{m^2}u(\widetilde t_n) + v(\widetilde t_n) -cu(\widetilde t_n),\quad v(\widetilde t_n)\in  \gamma^{-1}(u(\widetilde t_n)),$$ and we have~\eqref{parex01}, we arrive at
$$u(\widetilde t_n)\to  \gamma^+\quad\hbox{in } L^2(X,\nu),$$
which is impossible since the mass is preserved and we are taking $\displaystyle \gamma^-  < \frac{1}{\nu(X)} \int_X u_0 d\nu  <  \gamma^+ .$
Similarly we also arrive  to a contradiction if we suppose that $\lim \overline{\mu}(\widetilde t_n)=-\infty$.

Since $\{\overline{\mu}(\widetilde t_n)\}_n$ is bounded, then we have that there exist  a subsequence, that we denote equal, and a constant $\mu_\infty$ such that
$$\lim_n\overline{\mu}(\widetilde t_n)=\mu_\infty.$$
Therefore, from~\eqref{pere01}, we also have
$$\mu(\widetilde t_n,.)\to \mu_\infty\quad\hbox{in } L^2(X,\nu).$$
 And, from the convergences obtained, we easily arrive to
$$\mu_\infty   =  - \Delta_{m^2}u_\infty + v_\infty -cu_\infty, \quad v_\infty\in  \gamma^{-1}(u_\infty),$$
that is,
$u_\infty$ us a   stationary solution  of our Cahn-Hilliard problem. Hence we have proved~\eqref{gradlike}.
\end{proof}
\begin{remark}{\rm

\noindent (1) Let us point out that, for $\displaystyle F_3(r)=\frac{c}{2}(1 -r^2) + I_{[-1,1]}(r)$,   the assumption \eqref{CCOnd2} is natural,  otherwise the solution is trivial, see Remark~\ref{excon3}.

\noindent (2) Observe that $$\mathcal{E}(w)=\widetilde{\mathcal{E}}_\infty, \quad \hbox{for all} \ w \in \omega(u_0).$$
  In fact, given $w \in \omega(u_0)$, there exists $t_n \to + \infty$ such that $u(t_n) \to w$ in $L^2(X, \nu)$. Then, since $\mathcal{E}$ is lower semi-continuous, we have
$$\mathcal{E}(w) \leq \liminf_{n \to \infty} \mathcal{E}(u(t_n)) = \widetilde{\mathcal{E}}_\infty.$$
On the other hand, given $t >0$, if $t \leq t_n$, since $\mathcal{E}$ is non-inceasing, we have $$ \widetilde{\mathcal{E}}(t) \geq \widetilde{\mathcal{E}}(t_n) \geq \mathcal{E}(w).$$ Thus $\widetilde{\mathcal{E}}(t)$ is bounded from below on the orbit $\{ u(t) \ : \ t \geq 0 \}$,  and by using again the monotonicity of $\widetilde{\mathcal{E}}(t)$, we have   there exist $\lim_{t \to \infty}\widetilde{\mathcal{E}}(t) \geq \mathcal{E}(w)$. Therefore, $\mathcal{E}(w)=\widetilde{\mathcal{E}}_\infty$.
 $\blacksquare$
}
\end{remark}

For the obstacle potential $F_3$, it is interesting to know when,  for some $\nu$-measurable $D\subset X$, we have
$$\1_D-\1_{X\setminus D}\in \mathbb{E},$$
that is, when there exists equilibria solutions that divide the space in two pure phases without  interface between them. In the next result we will be that this happen under some geometrical condition on $D$. Consider here that we are dealing with
\begin{equation}\label{DNLS1penalty}
\left\{\begin{array}{lll} u_t(t,x) = \Delta_{m^1}  \mu(t,x) , \quad &(t,x) \in (0, \infty) \times X, \\[10pt] \mu(t,x)    \in  -\delta \Delta_{m^2}u(t,x) + \partial F(u(t,x)), \quad &(t,x) \in (0, \infty) \times X, \\[10pt] u(0,x) = u_0(x), \quad &x \in  X,\end{array} \right.
\end{equation}
with $\delta >0$ (for which the same existence and uniqueness result holds true).

\begin{proposition}\label{phases} Suppose we have Problem~\eqref{DNLS1penalty} with the obstacle potential $F_3$.  Let $D\subset X$ be  $\nu$-measurable such that $\displaystyle \int_X(\1_D-\1_{X\setminus D})d\nu=\int_Xu_0$. If \begin{equation}\label{tocom001}
1+\frac12\left(\sup_{x\in D}\mathcal{H}_{\partial D}^{m^2}(x)+ \sup_{x\in X\setminus D}\mathcal{H}_{\partial (X\setminus D)}^{m^2}(x)\right)
\le  \frac{c}{\delta},
\end{equation}
then
\begin{equation}\label{ole1}
\1_D-\1_{X\setminus D}\in \mathbb{E}.
\end{equation}
\end{proposition}

\begin{proof} To prove \eqref{ole1}, we need   the existence of a constant $\mu$ such
 $$ \mu+\delta\Delta_{m^2}(\1_D-\1_{X\setminus D})(x)+c(\1_D(x)-\1_{X\setminus D}(x))\in \gamma^{-1}(\1_D(x)-\1_{X\setminus D}(x)).$$
Therefore, since $\Delta_{m^2}(\1_D-\1_{X\setminus D})(x)=2m^2_x(D)-1   - (\1_D(x) - \1_{X \setminus D}(x)) $, we need the existence of a constant~$\mu$ such that
$$ \mu+\delta(2m^2_x(D)-1)+ (c-\delta)(\1_D(x)-\1_{X\setminus D}(x))\in \gamma^{-1}(\1_D(x)-\1_{X\setminus D}(x));$$
which is equivalent to ask for
$$\begin{array}{l}\mu+\delta(2m^2_x(D)-1)+  (c-\delta) \ge 0 \quad\hbox{if }x\in D,\\[10pt]
\mu+\delta(2m^2_x(D)-1)-  (c-\delta) \le 0 \quad\hbox{if }x\in X\setminus D.
\end{array}
$$
that is,
 $$ \begin{array}{l}\displaystyle\frac{\mu}{\delta}\ge  1-2m^2_x(D) -\left(\frac{c}{\delta}-1\right)  \quad\hbox{if }x\in D,\\[12pt]
\displaystyle \frac{\mu}{\delta}\le  1-2m^2_x(D) +\left(\frac{c}{\delta}-1\right)  \quad\hbox{if }x\in X\setminus D.
\end{array}
$$
Then, we can find a constant $\mu$ satisfying the above inequalities  if
   $$\displaystyle 1 + \sup_{x\in X\setminus D}m^2_x(D) - \inf_{x\in D}m^2_x(D)   \leq \frac{c}{\delta}.$$
Using that  $1-2m^2_x(D)=\mathcal{H}_{\partial D}^{m^2}(x)=-\mathcal{H}_{\partial (X\setminus D)}^{m^2}(x)$,  the above inequality is equivalent  to \eqref{tocom001}, and the proof concludes.
\end{proof}

\begin{remark}\rm   Observe that if $c \geq 2\delta$ then \eqref{tocom001} holds, and therefore
\begin{center}
\eqref{ole1} also holds for any $c \geq 2\delta$.
\end{center}
Then, an small $\delta$, or a large $c$ in the potential, ensures the existence of equilibria with only pure states regions.
 $\blacksquare$
\end{remark}

  In the next result  we are dealing with strong solutions of Problem~\eqref{DNLS1penalty}.

\begin{proposition} Assume we are under the assumptions of Theorem \ref{ASINTB}. If  \begin{equation}\label{fact001}c<\delta \,\hbox{gap}(-\Delta_{m^2}),
 \end{equation}
then
\begin{equation}\label{convergence1}  \omega(u_0)=\{\overline{u_0}\},
  \end{equation}
  where $\displaystyle\overline{u_0}=\frac{1}{\nu(X)}\int_Xu_0d\nu$, which is equivalent to have
$$\lim_{t \to +\infty} u(t) = \overline{u_0}.$$
\end{proposition}

\begin{proof}
Let $u_\infty, \widetilde u_\infty\in \omega(u_0)$. As in the proof of~Theorem~\ref{ASINTB}, we can find   constants $\mu$ and $\widetilde\mu_\infty$ such that
$$\mu_\infty   =  - \delta\Delta_{m^2}u_\infty + v_\infty -cu_\infty, \quad v_\infty\in  \gamma^{-1}(u_\infty),$$
and  $$ \widetilde \mu_\infty   =  - \delta\Delta_{m^2}\widetilde u_\infty + \widetilde v_\infty -c\widetilde u_\infty, \quad \widetilde v_\infty\in  \gamma^{-1}(\widetilde u_\infty).$$
Then,
$$ j_\gamma^*(\widetilde u_\infty)-j_\gamma^*(u_\infty)\ge \int_X(\mu_\infty+\delta\Delta_{m^2}u_\infty +cu_\infty)(\widetilde u_\infty-u_\infty) d\nu $$
and
$$j_\gamma^*(u_\infty)-j_\gamma^*(\widetilde  u_\infty)\ge \int_X(\widetilde  \mu_\infty+\delta\Delta_{m^2}\widetilde  u_\infty +c
\widetilde u_\infty)( u_\infty-\widetilde  u_\infty)d\nu.$$
Hence, since we have preservation of mass,
$$ j_\gamma^*(\widetilde u_\infty)-j_\gamma^*(u_\infty)\ge \int_X(\delta \Delta_{m^2}u_\infty +cu_\infty)(\widetilde u_\infty-u_\infty) d\nu$$
and
$$ \delta j_\gamma^*(u_\infty)-j_\gamma^*(\widetilde  u_\infty)\ge \int_X(\delta\Delta_{m^2}\widetilde  u_\infty +c
\widetilde u_\infty)( u_\infty-\widetilde  u_\infty) d\nu.$$
 Then adding both last expressions
  $$0\ge -\delta\int_{X}\Delta_m(\widetilde u_\infty-u_\infty)dm^2_x(y)(\widetilde u_\infty-u_\infty)d\nu(x)-c\int_X(\widetilde u_\infty-u_\infty)^2d\nu,$$
 and integrating by parts we get,
 $$ 0\ge \frac{\delta}{2}\int_{X\times X}|\nabla(\widetilde u_\infty-u_\infty)|^2d(\nu\otimes m^2_x)-c\int_X(\widetilde u_\infty-u_\infty)^2d\nu.$$
 Hence, since $\displaystyle\int_X(\widetilde u_\infty-u_\infty)d\nu =0$, from the Poincar\'{e} inequality, we obtain
 $$ \delta\hbox{gap}(-\Delta_{m^2})\int_X|\widetilde u_\infty-u_\infty|^2d\nu\le c\int_X|\widetilde u_\infty-u_\infty|^2d\nu.$$
 Therefore, since we are assuming \eqref{fact001}
  we have that $\widetilde u_\infty=u_\infty$,  that is, the omega limit set is a singleton,
 $$\omega(u_0)=\{u_\infty\}.$$ Now, since $\overline{u_0}\in ]\gamma^-,\gamma^+[$, there exists (a constant) $\beta\in\gamma^{-1}(\overline{u_0})$. Therefore we can take $[u_\infty, v_\infty,\mu_\infty]=[\overline{u_0},\beta,\beta-c\overline{u_0}]$ in the above computations, and we get~\eqref{convergence1}.
 \end{proof}

  Let us point out that by the above result we have that if \eqref{fact001} holds, then this model is not suitable for phase separation.

\subsection{The case of  finite weighted discrete graphs}

  Let $G = (V(G), E(G))$ be a finite graph and suppose that we have two sets of   weights $w^1_{xy}$, $w^2_{xy}$ such that
\begin{equation}\label{igualdegre}d_x:= \sum_{y\sim x} w^1_{xy} = \sum_{y\sim x} w^2_{xy}, \quad \hbox{for all} \ x \in V(G).
\end{equation}
So,
$$\nu_G(A):= \sum_{x \in A} d_x,  \quad A \subset V(G),$$
is a invariant measure for $(m^i_x)$, $i =1,2$, being
$$m^i_x:=  \frac{1}{d_x}\sum_{y \sim x} w^i_{xy}\,\delta_y, \quad i=1,2.$$
Then, we can consider the random walk spaces  $[V(G), m^i,\nu_G]$, $i =1,2$ (see Example \ref{example.graphs}), that we assume for both to be  $m^i$-connected.

We have that each random walk space  $[V(G), m^i,\nu_G]$, $i =1,2$, satisfies a Poincar\'{e} inequality, and also  \eqref{1258New} holds for any $u_0 \in L^2(V(G), \nu_G)$. Therefore,   by Theorem~\ref{ExistUniqStrong}, there exists a unique solution strong solution $u(t)$ of Problem \eqref{DNLS1} for any initial data $u_0 \in L^2(V(G), \nu_G)$.
Since $L^2(V(G), \nu_G)$ is finite dimensional,  if we assume that the potential   $\displaystyle j_\gamma^*(r)- \frac{c}{2}r^2 $ is bounded from bellow, then, by Corollary~\ref{ACOT11}, we have
$$\{u(t):t>0\} \hbox{ is bounded in $L^2(V(G), \nu_G)$,}$$
from where we get that
\begin{equation}\label{compact}\{u(t):t>0\} \hbox{ is relatively compact in $L^2(V(G), \nu_G)$.}
\end{equation}
And, therefore, the $\omega$-limit set
$\omega(u_0)$
is not empty.   Hence, as a consequence of Theorem~\ref{ASINTB}, we have the following result.

\begin{theorem}\label{ASINTBGraph} Let $G = (V(G), E(G))$ be a finite graph satisfying \eqref{igualdegre}.
Assume  that that the potential   $\displaystyle j_\gamma^*(r)- \frac{c}{2}r^2  $ is bounded from bellow.
Let $u_0 \in L^2(V, \nu_G)$ be and assume that
\begin{equation}\label{CCOnd2grap}  \gamma^-     < \frac{1}{\nu_G(V)} \int_V u_0 d\nu_G  <  \gamma^+ .
\end{equation}
Then,
\begin{equation}\label{gradlikegrap}
\omega(u_0) \subset \mathbb{E}.
\end{equation}

\end{theorem}

 \

\noindent {\bf Acknowledgments.} The authors have been partially supported  by the Spanish MCIU and FEDER, project PGC2018-094775-B-100  and by Conselleria d'Innovaci\'{o}, Universitats, Ci\`{e}ncia y Societat Digital, project AICO/2021/223.

 \end{document}